\documentclass{amsart}
\usepackage{array}
\usepackage{multirow}
\usepackage{caption}
\usepackage{subcaption}
\usepackage{graphicx}
\usepackage{epsfig}
\usepackage{amsfonts,amsmath}
\usepackage{amssymb}
\usepackage{latexsym}
\usepackage{epsf}
\usepackage{graphicx,color,graphics}
\usepackage{amsmath,amssymb}
\usepackage{dsfont}
\usepackage[latin1]{inputenc}  
\newtheorem{theorem}{Theorem}[section]

\newtheorem{lemma}[theorem]{Lemma}

\newtheorem{remark}{Remark}

\numberwithin{equation}{section}

\begin{document}

\title[\hfil Schr\"odinger equation with infinite memory]{Well-posedness and stability for Schr\"odinger
equations with infinite memory}

\pagenumbering{arabic}

\author[M. M. Cavalcanti]{M. M. Cavalcanti}
\author[V. N. Domingos Cavalcanti]{V. N. Domingos Cavalcanti}
\address[M. M. Cavalcanti,V. N. Domingos Cavalcanti]{Department of Mathematics, State University of Maring\'a, 
Maring\'a, PR 87020-900, Brazil}
\email{mmcavalcanti@uem.br, vncavalcanti@uem.br}

\author[A. Guesmia]{A. Guesmia}
\address[A. Guesmia]{Institut Elie Cartan de Lorraine, UMR 7502, Universit\'e de Lorraine\\
3 Rue Augustin Fresnel, BP 45112, 57073 Metz Cedex 03, France}
\address[A. Guesmia]{
Department of Mathematics and Statistics\\
King Fahd University of Petroleum and Minerals, Dhahran 31261, Saudi Arabia}
\email{aissa.guesmia@univ-lorraine.fr}

\author[M. Sep\'ulveda]{
M. Sep\'ulveda}
\address[M. Sep\'ulveda]{
DIM and CI$^2$MA, Universidad del Concepci\'{o}n, Concepci\'{o}n, Chile}
\thanks{Research of Mauricio Sepúlveda C. was supported FONDECYT grant no. 1180868, and by ANID-Chile through the project {\sc Centro de Modelamiento Matem\'atico} (AFB170001)  of the PIA Program: Concurso Apoyo a Centros Científicos y Tecnológicos de Excelencia con Financiamiento Basal.}
\email{mauricio@ing-mat.udec.cl}

\begin{abstract}
We study in this paper the well-posedness and stability for two linear Schr\"odinger equations in $d$-dimensional open bounded domain under Dirichlet boundary conditions with an infinite memory. First, we establish the well-posedness in the sens of semigroup theory. Then, a decay estimate depending on the smoothness of initial data and the arbitrarily growth at infinity of the relaxation function
is established for each equation with the help of multipliers method and some arguments devised in \cite{gues1} and \cite{gues2}.
\end{abstract}



\maketitle

\section{Introduction}

The subject of this paper is studying the existence and decay of solutions for the following two Schr\"odinger equations with infinite memory:  
\begin{align}\label{Eq1}
\begin{cases}
iy_{t} (x,t)+ a \Delta y(x,t)- i\displaystyle\int_0^{\infty}\,f (s) \Delta y (x,t-s) \,ds = 0, \quad &x\in \Omega,\,t\in \mathbb{R}_+^* :=(0,\infty),\\
y(x,t)=0,\quad &x\in \Gamma,\,t\in\mathbb{R}_+^* ,\\
y(x,-t)=y_0 (x,t),\quad &x\in \Omega,\,t\in \mathbb{R}_+ :=[0,\infty)
\end{cases}
\end{align}
and 
\begin{align}\label{Eq2}
\begin{cases}
iy_{t} (x,t)+ a \Delta y(x,t) + i\displaystyle\int_0^{\infty}\,f (s) y (x,t-s) \,ds = 0, \quad &x\in \Omega,\,t\in \mathbb{R}_+^* ,\\
y(x,t)=0,\quad &x\in \Gamma,\,t\in\mathbb{R}_+^* ,\\
y(x,-t)=y_0 (x,t),\quad &x\in \Omega,\,t\in \mathbb{R}_+ ,
\end{cases}
\end{align}
where the subscript $t$ denotes the derivative with respect to the time variable $t$, $\Delta$ is the laplacian operator with respect to the space variable $x$, $\Omega\subset\mathbb{R}^d$ is an open bounded domaine with a smooth boundary $\Gamma$, $d\in \mathbb{N}^{*}$, 
$a\in \mathbb{R}_+^{*}$, $f\,:\mathbb{R}_{+}\to \mathbb{R}$ is a given function, $y_0$ is a fixed initial data and $y$ is the unknown of \eqref{Eq1} and \eqref{Eq2}.
\vskip0,1truecm
We would like here to mention some known papers in connection with well-posedness and stability of Schr\"odinger type equations, which the subject of our paper. 
\vskip0,1truecm
When the infinite memory is replaced by a damping, equation \eqref{Eq1} in the presence or not in \eqref{Eq1}$_1$ of a semilinear term; that is 
\begin{align}\label{Eq01}
iy_{t} (x,t)+ a \Delta y(x,t)+b\vert y\vert^p y +icy = 0 
\end{align}
($a,\,p,\,c\in\mathbb{R}_{+}^*$ and $b\in\mathbb{R}$), has been widely studied in the literature, where it is known that Schr\"odinger equations are globally well-posed under some smallness conditions on $p$; see \cite{11}. In the particular case $a=1$, $p=2$, $b\in\{1 ,-1\}$ and the domain is bounded, the exponential stability of \eqref{Eq01} was proved in \cite{36} under some smoothness and smallness conditionns on the initial data. A generalization to the case of inhomogeneous Dirichlet boundary cnditions was given in \cite{27}, where the decay rate depends on the regularity of solutions. Some exact controllability results in both Dirichlet and Neumann boundary conditions cases are also known for \eqref{Eq01}, see \cite{31} ($p=2$ and $c=0$). For more general semilinearity: $p=2$ or not (with $a=1$, $b=c$ and the domain is unbounded), some global existence results of solutions as well as the bolw-up phenomena were obtained in \cite{26} for two sets of initial data.   
\vskip0,1truecm
In the cited papers above, a full damping was considered (that is $c\in\mathbb{R}_{+}^*$). The authors of \cite{8,9,10}  treated the case of locally distributed dampings; that is $c$ is a function on space variable and vanishes on some part of the domain. They proved that the expoential satbility holds true when $a=1$, $p=2$, $b\in\{1 ,-1\}$ and the domain is unbounded. In this case, and for some two dimensional domains, the controllability of the model was proved in \cite{21}. In the one dimensional unbounded damain case with $a=-1$, the authors of \cite{15} proved some stabilization and blow-up properties for \eqref{Eq01} depending on the nonlinearity power $p$.
\vskip0,1truecm
In \cite{ccct}, the authors studied the existence as well as the stability in $\mathbb{R}^d$ of \eqref{Eq01} with $a=1$, $b=-1$ and the damping coefficient $c$ is a function on both space and time variables and may vanish when time goes to infinity. Moreover, the uniqueness of solution is proved when $d\in\{1 ,2\}$. Similar results were obtained in \cite{1} and \cite{2} in $d$-dimensional Riemannian manifolds and nonlinear local damping $ic(x)g(y)$ (instead of $icy$) but with $b=0$, where $g$ is a given function satisfying some properties. The authors of \cite{07} considered in Riemannian manifolds two more general forms than \eqref{Eq01} by taking $f(\vert y\vert)$ and $(-\Delta)^k (c(x)y)$ instead of $b\vert y\vert^p$ and $cy$, respectively, where $f$ is a given function satisfying some properties and $k\in\{0,\frac{1}{2}\}$. They proved that, at infinity, the energy functional goes to zero if $k=0$ (weak dissipation), and converges exponentially to zero when $k=\frac{1}{2}$ (strong dissipation).
\vskip0,1truecm
There exist in the literature several well-posedness and stability (theoretical and numerical) results also for higher order
Schr\"odinger equations. In this direction, see, for example, \cite{007} and the references therein. 
\vskip0,1truecm
For other well-posedness, stability and blow-up results related to Schr\"odinger types equations cited above, we refer the readers to, for example, \cite{4}-\cite{7}, \cite{13}, \cite{14}, \cite{16}-\cite{20}, \cite{24}, \cite{28}-\cite{30}, \cite{32}, \cite{35}, \cite{37} and the references therein.
\vskip0,1truecm
Our goals in the present paper is studying the existence, uniqueness, regularity and decay of solutions for the two linear Schr\"odinger equations \eqref{Eq1} and \eqref{Eq2}, where the unique present dissipation is the one generated by the infinite memory
term. This situation is completely different from the ones considered in the literature and cited above, where the dissipation is generated by a (linear or nonlinear) damping. First, we establish the well-posedness (existence, uniqueness and smoothness of solutions) in the sens of semigroup theory. Then, a decay estimate depending on the smoothness of initial data and the arbitrarily growth at infinity of the relaxation function $f$ is established for each equation. These two decay estimates imply that any weak solution converges to zero at infinity. In the particular case where $-f^{\prime}$ converges exponentially to zero at infinity, our decay estimates lead to the decay rate $t^{-n}$, where $n\in\mathbb{N}^*$ depends on the regularity of initial data (see section 3). The proofs are based on the semigroup approach, the multipliers method and some arguments devised in \cite{gues1} and \cite{gues2}. 
\vskip0,1truecm
The paper is organized as follows. In section 2, we present our assumptions on the function $f$, state and prove the well-posedness of (\ref{Eq1}) and (\ref{Eq2}). In section 3, we consider some assumptions on the growth of $f$ at infinity, state and prove our stability results. We give some general comments in section 4.
Finally, in section 5, we give some numerical examples that graphically illustrate the theoretical results obtained.

\section{Preliminaries and well-posedness results}

In this section, we present and proof our well-posedness results for \eqref{Eq1} and \eqref{Eq2}. To simplify the formulations, the variables $x$, $t$ and $s$ are noted only when it is needed to avoid ambiguity. Let us use $\left\langle , \right\rangle$ and $\Vert \cdot\Vert$ to denote, respectively, the standard inner product in $L^2 (\Omega)$ and its generated norm given by
\begin{equation*}
\left\langle p, q\right\rangle =\int_{\Omega} p(x){\bar {q}}(x)\,dx\quad\hbox{and}\quad \Vert p\Vert = \left(\int_{\Omega} \vert p(x)\vert^2\,dx\right)^{\frac{1}{2}} .
\end{equation*}
In order to prove the well-posedness of \eqref{Eq1} and \eqref{Eq2} using the semigroup approach, and as in \cite{dafe}, we consider the varibale $\eta^t$ and its initial data $\eta^0$ given by
\begin{equation}\label{etaname}
\eta^t (x,s)= \displaystyle\int_{t-s}^t y(x,\tau)\,d\tau\quad\hbox{and}\quad\eta^0 (x, s)=\displaystyle\int_{0}^s y_0 (x,\tau)\,d\tau,\quad x\in \Omega,\,\, s,t\in \mathbb{R}_{+} .
\end{equation}
Direct computations show that the functional $\eta^t$ satisfies
\begin{align}\label{eta}
\begin{cases}
\eta^t_t (x,s)+\eta^t_s (x,s) =y(x,t),\quad & x\in \Omega,\,\, s,t\in\mathbb{R}_+^* ,\\
\eta^t (x,s) =0,\quad & x\in \Gamma,\,\, s,t\in\mathbb{R}_+^* ,\\
\eta^t (x,0) =0,\quad & x\in \Omega,\,\, t\in \mathbb{R}_{+} ,
\end{cases}
\end{align}
where the subscript $s$ denotes the derivative with respect to $s$. To express in term of $\eta^t$ the memory integrals in \eqref{Eq1} and \eqref{Eq2}, we assume the following hypothesis: 
\vskip0,1truecm
${\bf (H1)}$ Assume that the function $f$ is non-increasing such that 
\begin{equation}\label{Co1}
f\in C^2 (\mathbb{R}_{+}),\quad f(0)>0\quad\hbox{and}\quad \lim_{s\to\infty} f(s)=0.
\end{equation}
\vskip0,1truecm
We put $g=-f'$, so $g\in C^1 (\mathbb{R}_{+})$, $g$ is non-negative and 
\begin{equation*}
g_0 :=\displaystyle\int_0^{\infty} g(s)\,ds =f(0)\in \mathbb{R}_{+}^*.
\end{equation*}
On the other hand, by integrating with respect to $s$ and using \eqref{eta}$_3$ and the limit in \eqref{Co1}, we get 
\begin{equation*}
\displaystyle\int_0^{\infty}\,g (s) \Delta \eta^t \,ds=-\displaystyle\int_0^{\infty}\,f' (s) \Delta \eta^t \,ds
=\displaystyle\int_0^{\infty}\,f (s) \Delta \eta^t_s \,ds.
\end{equation*}
From the definition of $\eta^t$, we see that $\eta^t_s =y(t-s)$, consequently
\begin{equation*}
\displaystyle\int_0^{\infty}\,g (s) \Delta \eta^t \,ds=\displaystyle\int_0^{\infty}\,f (s) \Delta y(t-s) \,ds.
\end{equation*}
Similarly, we have
\begin{equation*}
\displaystyle\int_0^{\infty}\,g (s) \eta^t \,ds=\displaystyle\int_0^{\infty}\,f (s)y(t-s) \,ds.
\end{equation*}
Then the equations \eqref{Eq1}$_1$ and \eqref{Eq2}$_1$ can be rewritten, respectively, in the forms
\begin{equation}\label{Eq1+}
iy_{t} + a \Delta y - i\displaystyle\int_0^{\infty}\,g(s) \Delta \eta^t \,ds = 0\quad\hbox{and}\quad
iy_{t} + a \Delta y + i\displaystyle\int_0^{\infty}\,g(s) \eta^t \,ds = 0.
\end{equation}
We consider the variable $U$ and its initial data $U_0$ given by 
\begin{equation}\label{U}
U=(y,\eta^t )\quad\hbox{and}\quad U_0 =(y_0 (\cdot,0),\eta^0 ).
\end{equation}
Now, we can formulate the systems \eqref{Eq1} and \eqref{Eq2} in the following initial value problem:
\begin{align}\label{222}
\begin{cases}
U_t (t) = \mathcal{A}_j U(t),\quad t>0,\\
U(0) = U_{0},
\end{cases}
\end{align}
where $j=1$ in case \eqref{Eq1}, $j=2$ in case \eqref{Eq2} and the operators $\mathcal{A}_j$ are defined by
\begin{equation*}
\mathcal{A}_1 U =\left(
\begin{array}{c}
ia\Delta y +\displaystyle\int_0^{\infty}\,g (s) \Delta \eta^t \,ds\\
\\
y-\eta^t_s
\end{array}
\right)
\end{equation*}
and
\begin{equation*}
\mathcal{A}_2 U =\left(
\begin{array}{c}
ia\Delta y -\displaystyle\int_0^{\infty}\,g (s) \eta^t \,ds\\
\\
y-\eta^t_s
\end{array}
\right).
\end{equation*} 
\vskip0,1truecm
Let us consider the spaces
\begin{equation*}
L_{1} =\left\{v:\mathbb{R}_{+}\to H^1_0 (\Omega) ,\,\displaystyle\int_0^{\infty}\,g (s) \Vert \nabla v (s)\Vert^2 \,ds <\infty\right\}\,\,\hbox{and}\,\,L_{2} =\left\{v:\mathbb{R}_{+}\to L^2 (\Omega) ,\,\displaystyle\int_0^{\infty}\,g (s) \Vert v (s)\Vert^2 \,ds <\infty\right\},
\end{equation*}
equipped with the inner product
\begin{equation*}
\langle v ,w\rangle_{L_1} = \displaystyle\int_0^{\infty}\,g (s) \langle \nabla v (s) ,\nabla w (s)\rangle \,ds
\quad\hbox{and}\quad\langle v ,w\rangle_{L_2} = \displaystyle\int_0^{\infty}\,g (s) \langle v (s) ,w (s)\rangle \,ds,
\end{equation*}
and the energy space  
\begin{equation*}
\mathcal{H}_j =L^2 (\Omega)\times L_j ,\quad j=1,2
\end{equation*}
equipped with the inner product
\begin{equation*}
\langle (v_1 ,v_2 ) ,(w_1 ,w_2 )\rangle_{\mathcal{H}_j} = \langle v_1 ,w_1 \rangle +\langle v_2 ,w_2 \rangle_{L_j}.
\end{equation*}
The domain of $D(\mathcal{A}_j )$ is given by
\begin{equation*}
D(\mathcal{A}_j )=\left\{U\in \mathcal{H}_j ,\,\, \mathcal{A}_j U\in \mathcal{H}_j ,\,\,\eta^t (x,0)=0\right\},
\end{equation*}
more precisely,
\begin{equation*}
D(\mathcal{A}_1 )=\left\{U\in \mathcal{H}_1 ,\,\, y\in H^1_0 (\Omega),\,\,\eta^t_s\in L_1 ,\,\,\eta^t (x,0)=0,\,\,ia\Delta y+\displaystyle\int_0^{\infty}\,g (s)\Delta \eta^t\,ds\in L^2 (\Omega)\right\}
\end{equation*}
and 
\begin{equation*}
D(\mathcal{A}_2 )=\left\{U\in \mathcal{H}_2 ,\,\,\eta^t_s\in L_2 ,\,\,\eta^t (x,0)=0,\,\,\Delta y \in  L^2 (\Omega)\right\}.
\end{equation*} 
To get the well-posedness of \eqref{222}, we assume the following additional hypothesis: 
\vskip0,1truecm
${\bf (H2)}$ Assume that $g$ is non-increasing such that there exists a postive constant $\beta_0$ satisfying
\begin{equation}\label{2.230}
-\beta_0 g \leq g'.
\end{equation} 
\vskip0,1truecm
The well-posedness results for \eqref{222} are given in this theorem.
\vskip0,1truecm
\begin{theorem}\label{Theorem 1.1}
Assume that ${\bf (H1)}$ and ${\bf (H2)}$ hold. Then, for any $n\in \mathbb{N}$ and $U_0 \in D(\mathcal{A}_j^n )$, system \eqref{222} admits a unique solution $U$ satisfying
\begin{equation}
U\in \cap_{k=0}^n C^k \left(\mathbb{R}_{+} ;D(\mathcal{A}_j^{n-k} )\right). \label{exist1}
\end{equation}
\end{theorem}
\vskip0,1truecm
\begin{proof} We mention first that $\mathcal{H}_j$ is a Hilbert space and $\mathcal{A}_j$ is linear.
The proof of Theorem \ref{Theorem 1.1} relies then on the Lumer-Philips theorem by proving that the operator $\mathcal{A}_j$ is dissipative and $I-\mathcal{A}_j$ is surjective ($I$ denotes the identity operator); that is $-\mathcal{A}_j$ is maximal monotone. So $\mathcal{A}_j$ is the infinitesimal generator of a $C_0$ semigroup of contraction on $\mathcal{H}_j$ and its domain $D(\mathcal{A}_j )$ is dense in $\mathcal{H}_j$. The conclusion of Theorem \ref{Theorem 1.1} follows immediately (see \cite{liu00} and \cite{pazy}).
\vskip0,1truecm
Second, we prove that 
\begin{equation}
\Re\left\langle \mathcal{A}_1 U ,U \right\rangle_{\mathcal{H}_1 }=\frac{1}{2}\displaystyle\int_0^{\infty} g' (s)\left\Vert \nabla\eta^t\right\Vert^{2}\,ds\quad\hbox{and}\quad\Re\left\langle \mathcal{A}_2 U ,U \right\rangle_{\mathcal{H}_2 }=\frac{1}{2}\displaystyle\int_0^{\infty} g' (s)\left\Vert \eta^t\right\Vert^{2}\,ds, \label{dissp}
\end{equation}
where $\Re$ denotes the real part. Hence, $\mathcal{A}_j$ is dissipative, since $g$ is non-increasing and \eqref{2.230} guarantees the boundedness of the integrals in \eqref{dissp}. Using the definition of $\mathcal{A}_1$ and $\left\langle ,\right\rangle_{\mathcal{H}_1}$, integrating by parts and using the boundary condition, we get
\begin{equation}
\left\langle \mathcal{A}_1 U ,U \right\rangle_{\mathcal{H}_1 }=-ia\Vert \nabla y\Vert^2 
+\displaystyle\int_0^{\infty} g(s)\left(\left\langle \nabla y ,\nabla\eta^t \right\rangle-\left\langle \nabla\eta^t ,\nabla y \right\rangle-\left\langle \nabla\eta^t_s ,\nabla\eta^t \right\rangle\right)\,ds. \label{st5++}
\end{equation}
Direct computations imply that
\begin{equation*}
\left\langle \nabla\eta^t_s ,\nabla\eta^t\right\rangle = \frac{1}{2}\left(\Vert \nabla\eta^t \Vert^2 \right)_s + 
i\displaystyle\int_{\Omega}\left( \Re \nabla\eta^t \cdot\Im \nabla\eta^t_s -\Im \nabla\eta^t \cdot\Re \nabla\eta^t_s\right)\,dx
\end{equation*}
and 
\begin{equation*}
\left\langle \nabla y ,\nabla\eta^t \right\rangle-\left\langle \nabla\eta^t ,\nabla y \right\rangle =2i \Im \left\langle \nabla y ,\nabla\eta^t \right\rangle, 
\end{equation*}
where $\Im$ denotes the imaginary part. Exploiting these two equalities, we deduce from \eqref{st5++} that
\begin{equation}
\left\langle \mathcal{A}_1 U ,U \right\rangle_{\mathcal{H}_1 }=-ia\Vert \nabla y\Vert^2  
+2i \Im \displaystyle\int_0^{\infty}g(s)\left\langle \nabla y ,\nabla\eta^t \right\rangle \,ds\label{st5+}
\end{equation}
\begin{equation*}
+i\displaystyle\int_0^{\infty} g(s)\displaystyle\int_{\Omega}\left(\Im \nabla\eta^t \cdot\Re \nabla\eta^t_s -\Re \nabla\eta^t \cdot\Im \nabla\eta^t_s \right)\,dx\,ds -\frac{1}{2}\displaystyle\int_0^{\infty} g(s)\left(\Vert \nabla\eta^t \Vert^2 \right)_s\,ds .
\end{equation*}
Integrating the last integral in \eqref{st5+} with respect to $s$ and taking the real part of the obtained formula we get the first equality in \eqref{dissp}. The second equality in \eqref{dissp} can be obtained using exactely the same arguments, where we get instead of \eqref{st5+} 
\begin{equation*}
\left\langle \mathcal{A}_2 U ,U \right\rangle_{\mathcal{H}_2 }=-ia\Vert \nabla y\Vert^2 +2i \Im \displaystyle\int_0^{\infty} g(s)\left\langle y ,\eta^t \right\rangle\,ds 
\end{equation*}
\begin{equation*}
+i\displaystyle\int_0^{\infty} g(s)\displaystyle\int_{\Omega}\left(\Im \eta^t \Re \eta^t_s -\Re \eta^t \Im \eta^t_s \right)\,dx\,ds -\frac{1}{2}\displaystyle\int_0^{\infty} g(s)\left(\Vert\eta^t \Vert^2 \right)_s\,ds .
\end{equation*}  
\vskip0,1truecm
Third, we prove that $I-\mathcal{A}_j$ is surjective. Let 
$F=(f_1 ,f_2 ) \in \mathcal{H}_j$. We prove that there exists $U \in D\left( \mathcal{A}_j\right)$ satisfying
\begin{equation}
U-\mathcal{A}_j U=F.  \label{ZF}
\end{equation}
Let us consider the case $j=1$. The last equation in \eqref{ZF} is reduced to 
\begin{equation}
\eta^t_{s} + \eta^t = y + f_{2}. \label{z7*}
\end{equation}
Integrating with respect to $s$ and noting that $\eta^t$ should satisfy $\eta^t (x,0) = 0$, we get
\begin{equation}
\eta^t = (1 - e^{-\,s})y + \displaystyle\int_{0}^{s} e^{\tau - s}\,f_{2} (\tau)\ d\tau. \label{z7}
\end{equation}
The second equation in $\eqref{ZF}$ is reduced to
\begin{equation}\label{z5f5+}
y-ia\Delta y-\displaystyle\int_{0}^{\infty} g(s)\Delta \eta^t\,ds=f_1 .
\end{equation}
Multiplying \eqref{z5f5+} by ${\bar w}$, with $w\in H^1_0 (\Omega)$, integrating over $\Omega$ and using \eqref{z7}, we find
the variational formulation of $\eqref{z5f5+}$ given by
\begin{equation}\label{z5f5}
(g_1 +ia)\left\langle \nabla y ,\nabla w\right\rangle +\left\langle y ,w\right\rangle=\left\langle f_1 ,w\right\rangle -
\left\langle f_3 ,w\right\rangle_{L_1}, 
\end{equation}
where 
\begin{equation}\label{g1f3}
g_1 = \displaystyle\int_{0}^{\infty} (1-e^{-\,s} )\,g (s)\, ds\quad\hbox{and}\quad f_3 (s)= \displaystyle\int_{0}^{s} e^{\tau -\,s} \,f_2 (\tau)\, d\tau .
\end{equation}
We have, using the Fubini theorem and H\"older's inequality, we get
\begin{equation*}
\begin{array}{lll}
\displaystyle\int_{0}^{\infty}g (s)\left\Vert \nabla f_3\right\Vert^{2} \,ds
& = & \displaystyle\int_{0}^{\infty}g (s)\left\Vert \displaystyle\int_{0}^{s} e^{\tau -\,s} \,\nabla f_2 (\tau)\, d\tau \right\Vert^{2} \,ds \\
& \leq & \displaystyle\int_{0}^{\infty}e^{-\,2\,s}\,g (s)\left(\displaystyle\int_{0}^{s}e^{\tau}\ d\tau\right)\displaystyle\int_{0}^{s}e^{\tau}\,\Vert \nabla f_{2}(\tau)\Vert^{2}\, d\tau\, ds \\
\\
& \leq & \displaystyle\int_{0}^{\infty}e^{-\,s}\,(1 - e^{-\,s})\,g (s)\displaystyle\int_{0}^{s}e^{\tau}\,\Vert \nabla f_{2}(\tau)\Vert^{2}\,d\tau\,ds \\
\\
& \leq & \displaystyle\int_{0}^{\infty}e^{-\,s}\,g (s)\displaystyle\int_{0}^{s} e^{\tau}\,\Vert \nabla f_{2}(\tau)\Vert^{2}\, d\tau\, ds \\
\\
& \leq & \displaystyle\int_{0}^{\infty}e^{\tau}\,\Vert \nabla f_{2}(\tau)\Vert^{2} \displaystyle\int_{\tau}^{\infty}e^{-\,s}\,g (s) \,ds\, d\tau \\
\\
& \leq & \displaystyle\int_{0}^{\infty}e^{\tau}\,g (\tau)\,\Vert \nabla f_{2}(\tau)\Vert^{2} \displaystyle\int_{\tau}^{\infty}e^{-\,s}\,ds\,d\tau \\
\\
& \leq & \displaystyle\int_{0}^{\infty}g (\tau)\Vert \nabla f_{2}(\tau)\Vert^{2}\, d\tau =\Vert f_{2}\Vert_{L_1}^{2} < \infty,
\end{array}
\end{equation*} 
then $f_3 \in L_{1}$. Therefore, we see that, if $\eqref{z5f5+}$ admits a solution $y$ satisfying the required regularity in $D(\mathcal{A}_1 )$, then \eqref{z7} implies that $\eta$ exists and satisfies $\eta^t_{s} ,\,\eta^t\in L_{1}$. To prove the existence of $y$,
we notice that the form
\begin{equation*}
F_1 (v,w)=(g_1 +ia)\left\langle \nabla v ,\nabla w\right\rangle +\left\langle v ,w\right\rangle,\quad v,\,w\in H^1_0 (\Omega),
\end{equation*}
is bilinear, continuous and coercive, and the form
\begin{equation*}
F_2 (w)=\left\langle f_1 ,w\right\rangle -\left\langle f_3 ,w\right\rangle_{L_1}, \quad w\in H^1_0 (\Omega), 
\end{equation*}
is linear and continuous. For continuity of $F_1$ and $F_2$, we have just to apply the classical Poincar\'e's inequality: there exists 
$c_* >0$ such that
\begin{equation}\label{poincareine}
\Vert v\Vert^2 \leq c_* \Vert\nabla v\Vert^2 ,\quad v,\,w\in H^1_0 (\Omega).
\end{equation}
So using the Lax-Milgram theorem, we deduce that there exists a unique $y\in H^1_0 (\Omega)$ satisfying
\begin{equation*}
F_1 (y,w)=F_2 (w),\quad w\in H^1_0 (\Omega),
\end{equation*} 
which implies that \eqref{z5f5} holds. Hence, classical elliptic regularity arguments imply \eqref{z5f5+} and 
\begin{equation*}
ia\Delta y+\displaystyle\int_0^{\infty}\,g (s)\Delta \eta^t\,ds\in L^2 (\Omega).
\end{equation*}
This proves that $\eqref{ZF}$ in case $j=1$ has a unique solution $U\in D\left( \mathcal{A}_1\right)$.
The surjectivity of $I-\mathcal{A}_2$ can be proved in the same way, where in this case $\eta^t$ is defined in \eqref{z7}, the forms $F_1$ and $F_2$ are given by
\begin{equation*}
F_1 (v,w)=a\left\langle \nabla v ,\nabla w\right\rangle -i(1+g_1)\left\langle v ,w\right\rangle\quad\hbox{and}\quad F_2 (w)=-i\left\langle f_1 ,w\right\rangle +i\left\langle f_3 ,w\right\rangle_{L_2}, \quad v,\,w\in H^1_0 (\Omega)
\end{equation*}
and $g_1$ and $f_3$ are defined in \eqref{g1f3}. 
\end{proof}

\section{Stability results}

In this section, we present and prove our stability results for \eqref{222}, where the obtained decay estimate is valide for 
$U_0\in D(\mathcal{A}_1^{2n+2})$ in case \eqref{Eq1}, $U_0\in D(\mathcal{A}_2^{2n})$ in case \eqref{Eq2} and $n\in \mathbb{N}^*$.
We assume the following additional hypothesis on the growth of $g$ at infinity and the size of $y_0$: 
\vskip0,1truecm
${\bf (H3)}$ Assume that there exists a positive constant $\alpha_0$ and an increasing strictly convex function 
$G: \mathbb{R}_{+}\rightarrow  \mathbb{R}_{+}$ of class $C^{1}(\mathbb{R}_{+})\cap C^{2}(\mathbb{R}_{+}^*)$ satisfying 
\begin{equation}
\label{202}G(0) = G'(0) = 0\quad \mbox{and}\quad \lim_{t\rightarrow \infty}G'(t) = \infty
\end{equation}
such that
\begin{equation}\label{203+}
g^{\prime} \leq -\, \alpha_{0}\,g
\end{equation}
or
\begin{equation}\label{203}
\displaystyle\int_{0}^{\infty}\frac{s^2 g (s)}{G^{-1}\left(-\,g' (s)\right)}\ ds +\sup_{s\in\mathbb{R}_{+}}\frac{g (s)}{G^{-1}\left(-\,g' (s)\right)} < \infty .
\end{equation} 
Moreover, if \eqref{203+} does not hold, we assume that $y_0$ satisfies 
\begin{equation}\label{203*}
\sup_{t\in \mathbb{R}_{+}} \max_{k=0}^{n+1}\displaystyle\int_{t}^{\infty}\frac{g (s)}{G^{-1}\left(-\,g' (s)\right)}\left\Vert\displaystyle\int_{0}^{s-t}\nabla \partial_{s}^k y_0 (\cdot,\tau)\,d\tau\right\Vert^2 \, ds < \infty
\end{equation}
in case \eqref{Eq1}, where $\partial_{s}^k$ denotes the derivative of order $k$ with respect to $s$, and
\begin{equation}\label{203*+}
\sup_{t\in \mathbb{R}_{+}} \max_{k=0}^{n+1}\int_{t}^{\infty}\frac{g (s)}{G^{-1}\left(-\,g' (s)\right)}\left\Vert\displaystyle\int_{0}^{s-t} \partial_{s}^k y_0 (\cdot,\tau)\,d\tau\right\Vert^2\, ds < \infty
\end{equation}
in case \eqref{Eq2}.
\vskip0,1truecm
\begin{remark}\label{remark0}
1. Thanks to \eqref{203}, \eqref{203*} (resp. \eqref{203*+}) is valid if, for example, $\Vert\nabla \partial_{s}^k y_0\Vert^2$ (resp. $\Vert\partial_{s}^k y_0\Vert^2$), $k=0,1,\cdots,n+1$, are bounded with respect to $s$.
\vskip0,1truecm
2. The class of functions $g$ satisfying ${\bf (H1)}$, ${\bf (H2)}$ and ${\bf (H3)}$ is very wide and contains the ones which converge to zero exponentially like 
\begin{equation}\label{g2}
g_1 (s)=d_1 e^{-q_1 s}
\end{equation}
or at a slower rate like 
\begin{equation}\label{g1}
g_2 (s)=d_2 (1+s)^{-q_2}
\end{equation}
with $d_1,\,q_1 ,\,d_2 >0$ and $q_2 >3$. Conditions \eqref{2.230} and \eqref{203+} are satisfied by $g_1$ with $\beta_0 =\alpha_0 =q_1$, and conditions \eqref{2.230} and \eqref{203} are satisfied by $g_2$ with 
\begin{equation}\label{g1G}
\beta_0 =q_2 \quad\hbox{and}\quad G(s)=s^{p},\quad\hbox{for any}\,\, p >\dfrac{q_2 +1}{q_2 -3} . 
\end{equation}
\end{remark} 
\vskip0,1truecm
In order to announce our stability results, we consider the energy functionals $E_1$ and $E_2$ associated with \eqref{Eq1} and \eqref{Eq2}, respectively, and given by
\begin{equation}\label{En1}
E_1 (t) = \frac{1}{2} \left(\Vert y\Vert^2 +\displaystyle\int_0^{\infty} g(s)\Vert \nabla \eta^t\Vert^2 ds\right)
\end{equation}
and
\begin{equation}\label{En2}
E_2 (t) = \frac{1}{2} \left(\Vert y\Vert^2 +\displaystyle\int_0^{\infty} g(s)\Vert \eta^t\Vert^2 ds\right).
\end{equation} 
\\
\begin{theorem}\label{theorem1}
Assume that ${\bf (H1)}$, ${\bf (H2)}$ and ${\bf (H3)}$ hold. Let $n\in \mathbb{N}^*$, $U_0\in D(\mathcal{A}_1^{2n+2})$ in case \eqref{Eq1} and $U_0\in D(\mathcal{A}_2^{2n})$ in case \eqref{Eq2}. Then there exist positive constants $\alpha_{1,n}$ and $\alpha_{2,n}$ such that 
\begin{equation}\label{203**}
E_1 (t)\leq \alpha_{1,n} \,G_{n}\left(\frac{\alpha_{1,n}}{t}\right)\quad\hbox{and}\quad E_2 (t)\leq \alpha_{2,n} \,G_{n}\left(\frac{\alpha_{2,n}}{t}\right),\quad t\in \mathbb{R}_{+}^* ,
\end{equation}
where
\begin{equation}\label{203***} 
G_m (s)=G_1 (sG_{m-1}(s)),\,\, m=2,3,\cdots ,n,\,\, G_1 =G_0^{-1} \,\,\hbox{and}\,\,
G_{0} (s) =
\begin{cases}
s & \,\,\hbox{if \eqref{203+} holds}, \\
s\,G'(s) & \,\,\hbox{if \eqref{203} holds}.
\end{cases}
\end{equation} 
\end{theorem}
\vskip0,1truecm
\begin{remark}\label{remark2}
1. We see that $G_{n} (0)=0$, then \eqref{203**} implies that
\begin{equation}\label{limnablay} 
\lim_{t\to\infty} E_j (t) =0 . 
\end{equation}
By density of $D(\mathcal{A}_1^4 )$ and $D(\mathcal{A}_2^2 )$ in $\mathcal{H}_1$ and $\mathcal{H}_2$, respectively, we conclude that \eqref{limnablay} is valid, for any 
$U_0\in\mathcal{H}_j$.
\vskip0,1truecm
2. In case \eqref{203+}, $G_m (s)=s^{m}$, and so \eqref{203**} is reduced to, for some $\beta_{1,n} ,\,\beta_{2,n} >0$,
\begin{equation}\label{expon} 
E_1 (t)\leq \beta_{1,n} \,t^{-n}\quad\hbox{and}\quad E_2 (t)\leq \beta_{2,n} \,t^{-n}.
\end{equation} 
However, in case \eqref{203}, \eqref{203**} is weaker that \eqref{expon}. For the example \eqref{g1}, where \eqref{203} is satisfied with $G$ given in \eqref{g1G}, \eqref{203**} implies that there exist $\gamma_{1,n},\,\gamma_{2,n} >0$ such that 
\begin{equation*}
E_1 (t)\leq \gamma_{1,n} t^{-p_n} \quad\hbox{and}\quad E_2 (t)\leq \gamma_{2,n} t^{-p_n} ,
\end{equation*}
where $p_n =\sum_{m=1}^{n} p^{-m}$. Notice that $p_n\to n$ when $p\to 1$; that is when $q_2\to \infty$. This means that, if $g$ converges to zero at infinity faster than any polynomial, then the decay rate given in \eqref{203**} is arbitrarily close to $t^{-n}$. 
\end{remark} 
\vskip0,1truecm
We start the proof of \eqref{203**} by proving the following indentities for the derivatives of $E_1$ and $E_2$:
\vskip0,1truecm
\begin{lemma}\label{lemma2}
The energy functionals $E_1$ and $E_2$ satisfy
\begin{equation}
E_1^{\prime} (t)=\frac{1}{2} \int_0^{\infty} g^{\prime} (s)\Vert\nabla \eta^t\Vert^2 \,ds \label{St2}
\end{equation}
and
\begin{equation}
E_2^{\prime} (t)=\frac{1}{2} \int_0^{\infty} g^{\prime} (s)\Vert \eta^t\Vert^2 \,ds.\label{St3}
\end{equation}
\end{lemma} 
\vskip0,1truecm
\begin{proof}
We have
\begin{equation*}
E_j^{\prime} (t)=\Re\,\frac{1}{2} \left(\left\Vert U\right\Vert^2_{\mathcal{H}_j }\right)^{\prime} =\Re\,\left\langle U_t ,U \right\rangle_{\mathcal{H}_j } .
\end{equation*}
Then, using \eqref{222}$_1$ and \eqref{dissp}, \eqref{St2} and \eqref{St3} follow.
\end{proof}
\vskip0,1truecm
\begin{remark}\label{remark1}
1. Thanks to ${\bf (H2)}$, $E_1^{\prime}$ and $E_2^{\prime}$ are well-defined and non-positive, and so \eqref{222} is dissipative.
\vskip0,1truecm
2. Because $U_0\in D(\mathcal{A}_1^{2n+2})$ in case \eqref{Eq1} and $U_0\in D(\mathcal{A}_2^{2n})$ in case \eqref{Eq2} with $n\in \mathbb{N}^*$, we can define the following energy functionals of higher order $E_{j,k}$, $k=1,2,3,4$ if $j=1$, and $k=1,2$ if $j=2$: 
\begin{equation}\label{defejk}
E_{j,k} (t)=\frac{1}{2} \left\Vert \partial_t^k U\right\Vert^2_{\mathcal{H}_j} .
\end{equation}
Because \eqref{Eq1} and \eqref{Eq2} are linear and the coefficient $a$ does not depend on $t$, we get (as for \eqref{St2} and \eqref{St3}) 
\begin{equation}
E_{1,k}^{\prime} (t)=\frac{1}{2} \displaystyle\int_0^{\infty} g^{\prime} (s)\Vert\nabla \partial_t^k\eta_t\Vert^2 \,ds\quad\hbox{and}\quad E_{2,k}^{\prime} (t)=\frac{1}{2} \displaystyle\int_0^{\infty} g^{\prime} (s)\Vert\partial_t^k\eta^t\Vert^2 \,ds.\label{St2*}
\end{equation}
\end{remark} 
\vskip0,1truecm
\begin{lemma}\label{lemma3}
There exit positive constants $c_1$, $c_2$, ${\tilde c}_1$ and ${\tilde c}_2$  such that 
\begin{equation}
\Vert\nabla y\Vert^2\leq c_1\displaystyle\int_0^{\infty} g (s)\Vert\nabla\eta^t\Vert^2 \,ds +c_2 \displaystyle\int_{\Omega}\left(
\Im y\Re y_t -\Re y \Im y_t \right)\,dx\label{St2+}
\end{equation}
in case \eqref{Eq1}, and
\begin{equation}
\Vert\nabla y\Vert^2\leq {\tilde c}_1 \displaystyle\int_0^{\infty} g (s)\Vert\eta^t\Vert^2 \,ds +{\tilde c}_2 \displaystyle\int_{\Omega}\left(\Im y\Re y_t -\Re y \Im y_t \right)\,dx\label{St3+}
\end{equation}
in case \eqref{Eq2}.
\end{lemma} 
\vskip0,1truecm
\begin{proof}
Multiplying the first equation in \eqref{Eq1+} by ${\bar y}$, integrating over $\Omega$ and using the boundary condition, we get
\begin{equation}
a\Vert \nabla y \Vert^2 =i\left\langle y_t ,y \right\rangle +i\displaystyle\int_0^{\infty} g(s)\left\langle \nabla\eta^t ,\nabla y\right\rangle\,ds.\label{st7}
\end{equation}
Direct computations lead to
\begin{equation}
\left\langle y_t ,y \right\rangle=\frac{1}{2} \left(\Vert y\Vert^2 \right)_t +i\displaystyle\int_{\Omega}\left(
\Re y\Im y_t -\Im y\Re y_t \right)\,dx.\label{identyyt}
\end{equation}
On the other hand, applying H\"older's and Young's inequalities, we have, for any $\epsilon>0$,
\begin{equation}\label{valabsine}
\left\vert i\displaystyle\int_0^{\infty} g(s)\left\langle \nabla\eta^t ,\nabla y\right\rangle\,ds\right\vert\leq \epsilon \Vert \nabla y\Vert^2 +c_{\epsilon} \displaystyle\int_0^{\infty} g(s)\Vert\nabla\eta^t\Vert^2\,ds,
\end{equation}
where we denote by $c_{\epsilon}$ a positive constant depending on $\epsilon$. Combining \eqref{st7} and \eqref{identyyt}, taking the real part, using \eqref{valabsine} and choosing $\epsilon=\frac{a}{2}$, we deduce \eqref{St2+}.
\vskip0,1truecm
Similarly, multiplying the second equation in \eqref{Eq1+} by ${\bar y}$, integrating over $\Omega$ and using the boundary condition and \eqref{identyyt}, we find
\begin{equation}\label{secondinest}
a\Vert \nabla y \Vert^2 =i\displaystyle\int_0^{\infty} g(s)\left\langle \eta^t , y\right\rangle\,ds+\displaystyle\int_{\Omega}\left(
\Im y\Re y_t -\Re y\Im y_t\right)\,dx+\frac{i}{2} \left(\Vert y\Vert^2 \right)_t .
\end{equation}
Applying H\"older's, Young's and Poncar\'e's inequalities, we have, for any $\epsilon >0$,
\begin{equation}\label{secondinest+}
\left\vert i\displaystyle\int_0^{\infty} g(s)\left\langle \eta^t ,y\right\rangle\,ds\right\vert \leq \epsilon \Vert \nabla y\Vert^2 +c_{\epsilon} \displaystyle\int_0^{\infty} g(s)\Vert \eta^t\Vert^2\,ds.
\end{equation}
By taking the real part of \eqref{secondinest}, using \eqref{secondinest+} and choosing $\epsilon=\frac{a}{2}$, we obtain \eqref{St3+}.
\end{proof}
\vskip0,1truecm
Now, we prove the following estimations on the last integral in \eqref{St2+} and \eqref{St3+}:
\vskip0,1truecm
\begin{lemma}\label{Lemma4}
For any $\epsilon >0$, we have 
\begin{equation}\label{at}
\displaystyle\int_{\Omega}\left(\Im y \Re y_t -\Re y \Im y_t \right)\,dx\leq \epsilon \Vert\nabla y\Vert^2 +c_{\epsilon} \displaystyle\int_{0}^{\infty} g (s)\Vert\nabla \eta^t_{tt}\Vert^2 \,ds-c_{\epsilon} E_{1,1}^{\prime} (t)
\end{equation} 
in case \eqref{Eq1}, and
\begin{equation}\label{at+}
\displaystyle\int_{\Omega}\left(\Im y \Re y_t -\Re y \Im y_t \right)\,dx\leq \epsilon \Vert\nabla y\Vert^2 +c_{\epsilon} \displaystyle\int_{0}^{\infty} g (s)\Vert\eta^t_{tt}\Vert^2 \,ds -c_{\epsilon} E_{2,1}^{\prime} (t)
\end{equation}
in case \eqref{Eq2}. 
\end{lemma}
\vskip0,1truecm
\begin{proof}
We proceed as in \cite{fern1} for Timoshenko systems. Exploiting \eqref{eta} and integrating with respect to $s$, we have 
\begin{equation*}
\displaystyle\int_{\Omega}\left(\Im y \Re y_t -\Re y \Im y_t \right)\,dx= \frac{1}{g_0}\displaystyle\int_{\Omega} \Im y \displaystyle\int_0^{\infty} g(s)\Re\left(\eta^t_{tt}+\eta^t_{st}\right)\,ds\,dx -\frac{1}{g_0}\displaystyle\int_{\Omega} \Re y \displaystyle\int_0^{\infty} g(s)\Im\left(\eta^t_{tt}+\eta^t_{st}\right)\,ds\,dx. 
\end{equation*}
\begin{equation}\label{at0}
=  \frac{1}{g_0}\displaystyle\int_{\Omega} \Im y \displaystyle\int_0^{\infty} \Re\left(g(s)\eta^t_{tt} -g^{\prime} (s)\eta^t_{t}\right)\,ds\,dx-\frac{1}{g_0}\displaystyle\int_{\Omega} \Re y \displaystyle\int_0^{\infty} \Im\left(g(s)\eta^t_{tt} -g^{\prime} (s)\eta^t_{t}\right)\,ds\,dx. 
\end{equation} 
Using H\"older's, Young's and Poincar\'e's inequalities, we get, for any $\epsilon >0$, 
\begin{equation}\label{at1}
\displaystyle\int_{\Omega}\left(\Im y \Re y_t -\Re y \Im y_t \right)\,dx\leq \epsilon \Vert\nabla y\Vert^2 +c_{\epsilon}\displaystyle\int_0^{\infty}\left(g(s)\Vert\nabla \eta^t_{tt}\Vert^2 -g^{\prime} (s) \Vert\nabla \eta^t_{t}\Vert^2\right)\,ds.
\end{equation} 
Exploiting \eqref{St2*}$_1$ for $k=1$, we see that \eqref{at1} leads to \eqref{at}. Very similarly (using Poincar\'e's inequality only for $y$ and exploiting \eqref{St2*}$_2$ for $k=1$), \eqref{at+} is proved.
\end{proof}
\vskip0,1truecm
Now, choosing $\epsilon =\frac{1}{2c_2}$ in \eqref{at} and combining with \eqref{St2+}, we find, for some $c_3 >0$,
\begin{equation}\label{Stab1}
\Vert\nabla y\Vert^2\leq c_3 \displaystyle\int_{0}^{\infty} g (s)\Vert \nabla \eta^t \Vert^2 \,ds
+c_3 \displaystyle\int_{0}^{\infty} g (s)\Vert \nabla \eta^t_{tt}\Vert^2 \,ds-c_3 E_{1,1}^{\prime} (t) 
\end{equation} 
in case \eqref{Eq1}. And by combining \eqref{St3+} and \eqref{at+} with $\epsilon =\frac{1}{2{\tilde c}_2}$, we get, for some $c_4 >0$,
\begin{equation}\label{Stab2}
\Vert\nabla y\Vert^2\leq c_4 \displaystyle\int_{0}^{\infty} g (s)\Vert \eta^t\Vert^2 ds+c_4 \displaystyle\int_{0}^{\infty} g (s)\Vert \eta^t_{tt} \Vert^2 \,ds-c_4 E_{2,1}^{\prime} (t)
\end{equation}
in case \eqref{Eq2}. 
\vskip0,1truecm
\begin{remark}\label{remark31}
Using \eqref{2.230}, \eqref{defejk} and \eqref{St2*}$_1$ (for $k=1$), we conclude from \eqref{Stab1} that in case \eqref{Eq1}, for some 
$c_5 >0$, 
\begin{equation}\label{boundgrady}
\Vert\nabla y\Vert^2\leq c_5 (E_{1} (t)+E_{1,1} (t)+E_{1,2} (t))\leq c_5 (E_{1} (0)+E_{1,1} (0)+E_{1,2} (0)). 
\end{equation}
Therefore, using H\"older's inequality, we find, for $t\geq s\geq 0$, 
\begin{equation*}
\Vert \nabla\eta^t \Vert^2 =\left\Vert\displaystyle\int_{t-s}^{t}\nabla y(\cdot,\tau)\,d\tau\right\Vert^2 \leq s\displaystyle\int_{t-s}^{t}\Vert\nabla y(\cdot,\tau)\Vert^2\,d\tau \leq c_5 (E_{1} (0)+E_{1,1} (0)+E_{1,2} (0))s^2 .
\end{equation*}
For $s> t\geq 0$, using the same arguments, we have 
\begin{equation*}
\Vert \nabla\eta^t \Vert^2 =\left\Vert \displaystyle\int_{0}^{s-t} \nabla y_0 (\cdot,\tau)\,d\tau +\displaystyle\int_{0}^{t} \nabla y(\cdot,\tau)\,d\tau \right\Vert^2 
\end{equation*}
\begin{equation*}
\leq 2\left\Vert \displaystyle\int_{0}^{s-t} \nabla y_0 (\cdot,\tau)\,d\tau \right\Vert^2 +2\left\Vert\displaystyle\int_{0}^{t} \nabla y(\cdot,\tau)\,d\tau \right\Vert^2 
\end{equation*}
\begin{equation*}
\leq 2\left\Vert \displaystyle\int_{0}^{s-t} \nabla y_0 (\cdot,\tau)\,d\tau \right\Vert^2 +2c_5 (E_{1} (0)+E_{1,1} (0)+E_{1,2} (0)) s^2 .
\end{equation*}
Consequently
\begin{equation}\label{3.26}
\Vert \nabla\eta^t \Vert^2 \leq \left\{
\begin{array}{ll}
c_5 (E_{1} (0)+E_{1,1} (0)+E_{1,2} (0))s^2 & \,\,\hbox{if} \,\, 0\leq s\leq t, \\
\\
2\left\Vert \displaystyle\int_{0}^{s-t} \nabla y_0 (\cdot,\tau)\,d\tau \right\Vert^2 +2c_5 (E_{1} (0)+E_{1,1} (0)+E_{1,2} (0))s^2 &\,\,\hbox{if}\,\, s>t\geq 0\end{array}\right. :=M_1 (t,\,s). 
\end{equation}
In the case \eqref{Eq2}, and since $\Vert y\Vert^2$ is a part of $E_2$ and $E_2$ is non-increasing, we remark that, for $t\geq s\geq 0$, 
\begin{equation*}
\Vert \eta^t \Vert^2 =\left\Vert\displaystyle\int_{t-s}^{t} y(\cdot,\tau)\,d\tau\right\Vert^2 \leq s\displaystyle\int_{t-s}^{t}\Vert y(\cdot,\tau)\Vert^2\,d\tau \leq 2s\displaystyle\int_{t-s}^{t} E_2 (\tau)\,d\tau\leq 2E_2 (0)s^2 .  
\end{equation*}
For $s> t\geq 0$, we see that
\begin{equation*}
\Vert \eta^t \Vert^2 =\left\Vert \displaystyle\int_{0}^{s-t} y_0 (\cdot,\tau)\,d\tau +\displaystyle\int_{0}^{t} y(\cdot,\tau)\,d\tau \right\Vert^2 \leq 2\left\Vert \displaystyle\int_{0}^{s-t} y_0 (\cdot,\tau)\,d\tau \right\Vert^2 +4E_2 (0)s^2.
\end{equation*}
Hence
\begin{equation}\label{3.260}
\Vert \eta^t \Vert^2 \leq \left\{
\begin{array}{ll}
2E_2 (0)s^2 & \,\,\hbox{if} \,\, 0\leq s\leq t, \\
\\
2\left\Vert\displaystyle\int_{0}^{s-t} y_0 (\cdot,\tau)\,d\tau \right\Vert^2 +4E_2 (0)s^2 &\,\,\hbox{if}\,\, s>t\geq 0
\end{array}\right. :=M_2 (t,\,s).
\end{equation} 
Similarly to \eqref{3.26} and \eqref{3.260} and since $E_{1,2}$ and $E_{2,2}$ are non-increasing, we have, for some $c_6 >0$,  
\begin{equation}\label{m1tilde} 
\Vert \nabla\eta^t_{tt} \Vert^2 \leq \left\{
\begin{array}{ll}
c_6 (E_{1,2} (0)+E_{1,3} (0)+E_{1,4} (0)) s^2 & \quad\hbox{if} \quad 0\leq s\leq t, \\
\\
2\left\Vert \displaystyle\int_{0}^{s-t} \nabla \partial_{s}^2 y_0 (\cdot,\tau)\,d\tau \right\Vert^2 +2c_6 (E_{1,2} (0)+E_{1,3} (0)+E_{1,4} (0)) s^2 &\,\,\hbox{if}\,\, s>t\geq 0
\end{array}\right. :={\tilde M}_1 (t,\,s) 
\end{equation}
in the case \eqref{Eq1}, and
\begin{equation}\label{m2tilde}
\Vert \eta^t_{tt} \Vert^2 \leq \left\{
\begin{array}{ll}
2E_{2,2} (0)s^2 & \quad\hbox{if} \quad 0\leq s\leq t, \\
\\
2\left\Vert \displaystyle\int_{0}^{s-t}\partial_{s}^2 y_0 (\cdot,\tau)\,d\tau \right\Vert^2 +4E_{2,2} (0)s^2 &\quad\hbox{if}\quad s>t\geq 0
\end{array}\right. :={\tilde M}_2 (t,\,s)
\end{equation}
in the case \eqref{Eq2}. The inequalities \eqref{3.26}, \eqref{3.260}, \eqref{m1tilde} and \eqref{m2tilde} will be used in the proof of the next lemma in order to estimate the integrals in \eqref{Stab1} and \eqref{Stab2}. This lemma was introduced in \cite{gues1} and improved in \cite{gues2}. Notice that we have used energies of higher order up to $E_{1,4}$ in case \eqref{Eq1}, and up to $E_{2,2}$ in case \eqref{Eq2}; this why we need initial data $U_0\in D(\mathcal{A}_1^{2n+2})$ in case \eqref{Eq1} and $U_0\in D(\mathcal{A}_2^{2n})$ in case \eqref{Eq2} with $n\in \mathbb{N}^*$. 
\end{remark}  
\vskip0,1truecm
\begin{lemma}\label{lemma5} 
There exist positive constants $d_{1}$, ${\tilde d}_{1}$, $d_{2}$ and ${\tilde d}_{2}$ such that, for any $\varepsilon_{0}>0$, the following inequalities hold:
\begin{equation}
\frac{G_{0}\,(\varepsilon_{0}E_1 (t))}{\varepsilon_{0}E_1 (t)}\displaystyle\int_{0}^{\infty}g (s)\Vert \nabla \eta^t\Vert^2 \,ds\leq -\,d_{1}\,E_1^{\prime}(t)+ d_{1}\,G_{0}\,(\varepsilon_{0}E_1 (t)),\label{3.23}
\end{equation}
\begin{equation}
\frac{G_{0}\,(\varepsilon_{0}E_1 (t))}{\varepsilon_{0}E_1 (t)}\displaystyle\int_{0}^{\infty}g (s)\Vert \nabla \eta^t_{tt}\Vert^2 \,ds\leq -\,{\tilde d}_1\,E_{1,2}^{\prime}(t)+ {\tilde d}_1\,G_{0}\,(\varepsilon_{0}E_1 (t)),\label{3.23*}
\end{equation}
\begin{equation}
\frac{G_{0}\,(\varepsilon_{0}E_2 (t))}{\varepsilon_{0}E_2 (t)}\displaystyle\int_{0}^{\infty}g (s)\Vert \eta^t\Vert^2 \,ds\leq -\,d_{2}\,E_2^{\prime}(t)+ d_{2}\,G_{0}\,(\varepsilon_{0}E_2 (t))\label{3.23+}
\end{equation}
and
\begin{equation}
\frac{G_{0}\,(\varepsilon_{0}E_2 (t))}{\varepsilon_{0}E_2 (t)}\displaystyle\int_{0}^{\infty}g (s)\Vert \eta^t_{tt}\Vert^2 \,ds\leq -\,{\tilde d}_{2}\,E_{2,2}^{\prime}(t)+ {\tilde d}_{2}\,G_{0}\,(\varepsilon_{0}E_2 (t)),\label{3.23+*}
\end{equation}
where $G_0$ is defined in \eqref{203***}.
\end{lemma}
\vskip0,1truecm
\begin{proof} 
If \eqref{203+} holds, then \eqref{St2}, \eqref{St3} and \eqref{St2*} (for $k=2$) lead to
\begin{equation}\label{408}
\displaystyle\int_{0}^{\infty}g (s)\Vert \nabla \eta^t\Vert^2 \, ds\leq -\, \frac{2}{\alpha_{0}} E_1^{\prime}(t),\quad \displaystyle\int_{0}^{\infty}g (s)\Vert \nabla \eta^t_{tt}\Vert^2 \, ds\leq -\, \frac{2}{\alpha_{0}} E_{1,2}^{\prime}(t),
\end{equation}
\begin{equation}\label{408+}
\displaystyle\int_{0}^{\infty}g (s)\Vert \eta^t\Vert^2 \,ds\leq -\, \frac{2}{\alpha_{0}} E_2^{\prime}(t)\quad\hbox{and}\quad \displaystyle\int_{0}^{\infty}g (s)\Vert \eta^t_{tt}\Vert^2 \,ds\leq -\, \frac{2}{\alpha_{0}} E_{2,2}^{\prime}(t).
\end{equation}
So \eqref{3.23}-\eqref{3.23+*} hold with $d_{1} = {\tilde d}_{1} =d_2 ={\tilde d}_{2} = \frac{2}{\alpha_{0}}$ and $G_{0}(s)=s$.
\vskip0,1truecm
When \eqref{203+} does not hold and \eqref{203} is satisfied, we note first that, without loss of generality, we can assume that 
$E_1 >0$, $E_2 >0$ and $g^{\prime}<0$ on $\mathbb{R}_+$. Otherwise, if $E_1 (t_{1})=0$ and $E_2 (t_{2})=0$, for at least $t_{1},\,t_2\in \mathbb{R}_+$, then $E_1 (t)=0$, for all $t\geq t_{1}$, and $E_2 (t)=0$, for all $t\geq t_{2}$, since $E_1$ and $E_2$ are non-negative and non-increasing, and consequently, \eqref{203**} is satisfied, since $E_1$ and $E_2$ are bounded. And if $g^{\prime}<0$ is not satisfied on
$\mathbb{R}_+$, then there exists $s_0\in \mathbb{R}_+$ such that $g^{\prime} (s_0)=0$ and $g^{\prime} <0$ on $(0,s_0)$, since 
$g^{\prime}\in C(\mathbb{R}_+)$. Therefore \eqref{203} implies that $g(s_0)=0$, and so $g(s)=0$, for all $s\geq s_0$, since $g$ is non-negative and non-increasing. Consequently, the integrals on $\mathbb{R}_+$ in \eqref{Stab1} and \eqref{Stab2} are reduced to integrals on $(0,s_0)$ and $g^{\prime} <0$ on $(0,s_0)$.
\vskip0,1truecm
Let $\tau_{1}(t,\,s),\,\tau_{2}(t,\,s),\,\theta_{1} (t,\,s),\,\theta_{2} (t,\,s)>0$, $\varepsilon_{0} >0$ (which will be fixed later on) and $K(s)={\frac{s}{{G^{-\,1}(s)}}}$, for $s>0$. The hypothesis ${\bf (H3)}$ implies that 
\begin{equation*}
\lim_{s\to 0^{+}}\dfrac{s}{G^{-1}(s)}=\lim_{\tau\to 0^{+}}\dfrac{G(\tau)}{\tau}=G'(0)=0,
\end{equation*}
then $K(0)=0$. The function $K$ is non-decreasing. Indeed, the fact that $G^{-1}$ is concave and $G^{-1} (0)=0$ implies that, for any $0\leq s_{1} < s_{2}$,
\begin{equation*}
K(s_{1}) = \dfrac{s_1}{G^{-\,1}\left(\frac{s_{1}}{s_{2}}s_{2} + \left(1 - \dfrac{s_{1}}{s_{2}}\right)0\right)} 
\leq \dfrac{s_{1}}{\dfrac{s_{1}}{s_{2}}\,G^{-1}(s_{2}) + \left(1 - \dfrac{s_{1}}{s_{2}}\right)G^{-1} (0)} = \dfrac{s_{2}}{G^{-1}(s_{2})} = K(s_{2} ).
\end{equation*} 
Then, using \eqref{3.26} and \eqref{3.260},
\begin{equation}
K\left(-\,\theta_{1} (t,\,s)\,g^{\prime}(s) \Vert\nabla \eta^t\Vert^2\right) \leq K\left(-\,M_1 (t,\,s) \,\theta_{1}(t,\,s)\,g^{\prime}(s)\right) \label{3.2601} 
\end{equation}
and
\begin{equation}
K\left(-\,\theta_{2} (t,\,s)\,g^{\prime}(s) \Vert \eta^t\Vert^2\right) \leq K\left(-\,M_2 (t,\,s) \,\theta_{2}(t,\,s)\,g^{\prime}(s)\right). \label{3.2602} 
\end{equation}
Using \eqref{3.2601}, we arrive at 
\begin{equation*}
\int_{0}^{\infty}g (s) \Vert\nabla \eta^t\Vert^2 ds = \dfrac{1}{G'(\varepsilon_{0}\,E_1 (t))}\int_{0}^{\infty} \dfrac{1}{\tau_{1}(t,\,s)}G^{-\,1}\left(-\,\theta_{1} (t,\,s)\,g'(s)\Vert\nabla \eta^t\Vert^2\right) 
\end{equation*} 
\begin{equation*}
\times\dfrac{\tau_{1}(t,\,s)\,G'(\varepsilon_{0}\,E_1 (t))\,g(s)}{-\,\theta_{1} (t,\,s)\,g'(s)}\ K\left(-\,\theta_{1} (t,\,s)\,g' (s)
\Vert\nabla \eta^t\Vert^2\right)ds
\end{equation*} 
\begin{equation*}   
\leq \dfrac{1}{G'(\varepsilon_{0}\,E_1 (t))}\int_{0}^{\infty}\dfrac{1}{\tau_{1}(t,\,s)}\,G^{-\,1}\left(-\,\theta_{1} (t,\,s)\,g' (s)\Vert\nabla \eta^t\Vert^2\right)
\end{equation*} 
\begin{equation*}  
\times\dfrac{\tau_{1}(t,\,s)G'(\varepsilon_{0}\,E_1 (t))\,g (s)}{-\,\theta_{1} (t,\,s)\,g' (s)}\ K\left(-\,M_1 (t,\,s)\,\theta_{1} (t,\,s)\,g' (s)\right)ds
\end{equation*} 
\begin{equation*}
\leq\dfrac{1}{G' (\varepsilon_{0}\,E_1 (t))}\int_{0}^{\infty}\dfrac{1}{\tau_{1}\,(t,\,s)}\,G^{-1}\left(-\,\theta_{1} (t,\,s)\,g' (s)\Vert\nabla \eta^t\Vert^2\right)
\end{equation*} 
\begin{equation*}
\times \frac{M_1 (t,\,s)\,\tau_{1}(t,\,s)\,G'(\varepsilon_{0} \,E_1 (t))\,g (s)}{G^{-\,1}\,(-\,M_1 (t,\,s)\,\theta_{1} (t,\,s)\,g' (s))}\ ds.
\end{equation*}
Let $G^{*} (s)=\sup_{\tau\in \mathbb{R}_{+}} \{s\,\tau - G (\tau)\}$, for $s\in \mathbb{R}_{+}$, denote the dual function of $G$. From the hypothesis ${\bf (H3)}$, we see that
\begin{equation*}
G^{*}(s) = s\,(G')^{-1} (s) - G ((G')^{-1} (s)),\quad s\in \mathbb{R}_{+} .
\end{equation*}
Using Young's inequality: $s_{1}\,s_{2}\leq G(s_{1}) + G^{*}(s_{2})$, for 
\begin{equation*}
s_{1} = G^{-1}\left(-\,\theta_{1} (t,\,s)\,g^{\prime}(s)\Vert\nabla \eta^t\Vert^2\right)\quad\mbox{and}\quad s_{2} ={\dfrac{{M_1 (t,\,s)\,\tau_{1}(t,\,s)\,G^{\prime}(\varepsilon_{0}\,E_1 (t))\,g (s)}}{{G^{-\,1}(-\,M_1 (t,\,s)\,\theta_{1} (t,\,s) g^{\prime }(s))}}}, 
\end{equation*}
we obtain
\begin{equation*}
\int_{0}^{\infty}g (s)\Vert\nabla \eta^t\Vert^2\, ds \leq  \dfrac{1}{G^{\prime}(\varepsilon_{0} \,E_1 (t))}\int_{0}^{\infty}\dfrac{{-\,\theta_{1} (t,\,s)}}{{\tau_{1}(t,\,s)}}g^{\prime}(s)\Vert\nabla \eta^t\Vert^2\, ds 
\end{equation*} 
\begin{equation*}
+ \dfrac{1}{G^{\prime}(\varepsilon_{0}\,E_1 (t))}\int_{0}^{\infty}\dfrac{1}{\tau_{1}(t,\,s)}\ G^{*}\left({\frac{{M_1 (t,\,s)\,\tau_{1}(t,\,s)\,G^{\prime}(\varepsilon_{0}\,E_1 (t))\,g (s)}}{{G^{-1}(-\,M_1 (t,\,s)\,\theta_{1} (t,\,s)\,g^{\prime} (s))}}}\right)ds. 
\end{equation*}
Using the fact that $G^{*}(s)\leq s\,(G')^{-\,1} (s)$, we get 
\begin{equation*}
\int_{0}^{\infty}g (s)\Vert\nabla \eta^t\Vert^2\, ds  \leq  \dfrac{-\,1}{G^{\prime }(\varepsilon_{0}\,E_1 (t))}\int_{0}^{\infty}\dfrac{{\theta_{1} (t,\,s)}}{{\tau_{1}(t,\,s)}} g^{\prime }(s)\Vert\nabla \eta^t\Vert^2\, ds 
\end{equation*}
\begin{equation*}
\mbox{} + \int_{0}^{\infty}{\frac{{M_1 (t,\,s)\,g (s)}}{{G^{-\,1}(-\,M_1 (t,\,s)\,\theta_{1} (t,\,s)\,g^{\prime}(s))}}}\,(G')^{-\,1}\left({\frac{{M_1 (t,\,s)\,\tau_{1}(t,\,s)\,G^{\prime}(\varepsilon_{0}\,E_1 (t))\,g (s)}}{{G^{-\,1}(-\,M_1 (t,\,s)\,\theta_{1} (t,\,s)\,g^{\prime}(s))}}}\right)ds. 
\end{equation*}
Then, using the fact that $(G')^{-1}$ is non-decreasing and choosing $\theta_{1} (t,\,s) = {\frac{{1}}{{M_1 (t,\,s)}}}$, we find 
\begin{equation*}
\int_{0}^{\infty}g (s)\Vert\nabla \eta^t\Vert^2\, ds \leq \dfrac{-\,1}{G^{\prime}(\varepsilon_{0}\,E_1 (t))}\int_{0}^{\infty} \dfrac{1}{M_1 (t,\,s)\,\tau_{1}(t,\,s)}\,g^{\prime}(s)\Vert\nabla \eta^t\Vert^2\, ds  
\end{equation*} 
\begin{equation*}
+ \int_{0}^{\infty}{\dfrac{{M_1 (t,\,s)\,g (s)}}{{G^{-\,1}(-\,g^{\prime}(s))}}}\ (G')^{-1}\left(m_{0}\,M_1 (t,\,s)\,\tau_{1}(t,\,s)\,G^{\prime}(\varepsilon_{0}\,E_1 (t))\right)ds,  
\end{equation*}
where $m_{0} = \sup_{s\in \mathbb{R}_{+}}{\frac{{g (s)}}{{G^{-1}(-\,g^{\prime} (s))}}}$ ($m_{0}$ exists according to \eqref{203}). Due to \eqref{203} and the restriction \eqref{203*} on $y_{0}$ (for $k=0$), we have
\begin{equation*}
\sup_{t\in \mathbb{R}_{+}}\displaystyle\int_0^{\infty}{\frac{{M_1 (t,\,s)\,g (s)}}{{G^{-\,1}(-\,g^{\prime}(s))}}}\, 
ds=:m_{1} <\infty .
\end{equation*} 
Therefore, choosing $\tau_{1}(t,\,s)= {\frac{{1}}{{m_{0}\,M_1 (t,\,s)}}}$ and using \eqref{St2}, we obtain 
\begin{equation}\label{00000}
\int_{0}^{\infty}g (s)\Vert\nabla \eta^t\Vert^2\, ds  
\leq \dfrac{{-\,m_{0}}}{{G^{\prime}(\varepsilon_{0}\,E_1 (t))}}\int_{0}^{\infty}g^{\prime}(s)\Vert\nabla \eta^t\Vert^2\, ds + \varepsilon_{0}\,E_1 (t)\int_{0}^{\infty}{\frac{{M_1 (t,\,s)\,g (s)}}{{G^{-\,1}(-\,g^{\prime}(s))}}}\, ds 
\end{equation} 
\begin{equation*}
\leq \dfrac{{-\,2\,m_{0}}}{{G^{\prime}(\varepsilon_{0}\,E_1 (t))}}\,E_1^{\prime}(t) + m_{1}\,\varepsilon_{0}\,E_1 (t),
\end{equation*}
which, by multiplying \eqref{00000} by $G^{\prime}(\varepsilon_{0}\,E_1 (t))=\frac{G_0 (\varepsilon_{0}\,E_1 (t))}{\varepsilon_{0}\,E_1 (t)}$, gives \eqref{3.23} with $d_{1} = \max \{2\,m_{0},\,m_{1}\}$. Repeating the same arguments with $E_2$, $\Vert \eta^t\Vert^2$, $\tau_2$ and $\theta_2$ instead of $E_1$, $\Vert\nabla \eta^t\Vert^2$, $\tau_1$ and $\theta_1$, respectively, and using \eqref{203*+} (for $k=0$), \eqref{St3} and \eqref{3.2602}, we get \eqref{3.23+} with $d_{2} = \max \{2\,m_{0},\,m_{2}\}$, where  
\begin{equation*}
m_{2} =\sup_{t\in \mathbb{R}_{+}}\displaystyle\int_0^{\infty}{\frac{{M_2 (t,\,s)\,g (s)}}{{G^{-\,1}(-\,g^{\prime}(s))}}}\, ds,\quad\tau_{2}(t,\,s)= {\frac{{1}}{{m_{0}\,M_2 (t,\,s)}}}\quad\hbox{and}\quad \theta_{2} (t,\,s)= {\frac{{1}}{{M_2 (t,\,s)}}}.
\end{equation*}
As for \eqref{3.2601} and \eqref{3.2602},  
\begin{equation*}
K\left(-\,{\tilde\theta}_{1} (t,\,s)\,g^{\prime}(s) \Vert\eta^t_{tt}\Vert^2\right) \leq K\left(-\,{\tilde M}_1 (t,\,s) \,{\tilde\theta}_{1}(t,\,s)\,g^{\prime}(s)\right) 
\end{equation*}
and
\begin{equation*}
K\left(-\,{\tilde\theta}_{2} (t,\,s)\,g^{\prime}(s) \Vert \eta^t_{tt}\Vert^2\right) \leq K\left(-\,{\tilde M}_2 (t,\,s) \,{\tilde\theta}_{2}(t,\,s)\,g^{\prime}(s)\right), 
\end{equation*}
for any positive functions ${\tilde\theta}_{1}$ and ${\tilde\theta}_{2}$, where ${\tilde M}_1$ and ${\tilde M}_2$ are defined in \eqref{m1tilde} and \eqref{m2tilde}. Consequently, using the above two inequalities and arguing as for \eqref{00000} with ${\tilde\tau}_{j}$ and ${\tilde\theta}_{j}$ instead of $\tau_{j}$ and $\theta_{j}$, respectively, $j=1,2$, we deduce \eqref{3.23*} and \eqref{3.23+*} with 
${\tilde d}_j = \max \{2\,m_{0},\,{\tilde m}_{j}\}$, where  
\begin{equation*}
{\tilde m}_{j} =\sup_{t\in \mathbb{R}_{+}}\displaystyle\int_0^{\infty}{\frac{{{\tilde M}_j (t,\,s)\,g (s)}}{{G^{-\,1}(-\,g^{\prime}(s))}}}\, ds,\quad {\tilde\tau}_{j} (t,\,s)= {\frac{{1}}{{m_{0}\,{\tilde M}_j (t,\,s)}}}\quad\hbox{and}\quad {\tilde\theta}_{j} (t,\,s)= {\frac{{1}}{{{\tilde M}_j (t,\,s)}}}.
\end{equation*}
\end{proof}
\vskip0,1truecm
Now, using \eqref{poincareine} and the definition of $E_1$, we see that
\begin{equation}
\frac{2}{c_*}E_1 (t)\leq \Vert \nabla y\Vert^2 +\frac{1}{c_*}\displaystyle\int_{0}^{\infty}g (s)\Vert \nabla \eta^t\Vert^2 \,ds, \label{finine1}
\end{equation}
therefore, multiplying \eqref{finine1} by $\frac{G_{0}\,(\varepsilon_{0}E_1 (t))}{\varepsilon_{0}E_1 (t)}$ and using \eqref{Stab1}, we find 
\begin{equation}
\frac{2}{\epsilon_0 c_*}G_{0}\,(\varepsilon_{0}E_1 (t))\leq \left(\frac{1}{c_*} +c_3\right)\frac{G_{0}\,(\varepsilon_{0}E_1 (t))}{\varepsilon_{0}E_1 (t)}\displaystyle\int_{0}^{\infty}g (s)\left(\Vert \nabla \eta^t\Vert^2 +\Vert \nabla \eta^t_{tt}\Vert^2\right)\,ds
 \label{finine2}
\end{equation}
\begin{equation*}
-c_3 \frac{G_{0}\,(\varepsilon_{0}E_1 (t))}{\varepsilon_{0} E_1 (t)} E_{1,1}^{\prime} (t),
\end{equation*}
then, combining \eqref{finine2} with \eqref{3.23} and \eqref{3.23*}, we get 	
\begin{equation}\label{Stab11} 
\left[\frac{2}{c_*}-\epsilon_0 \left( c_3 +\frac{1}{c_*}\right) \left(d_1 +{\tilde d}_1\right)\right] G_{0}\,(\varepsilon_{0}E_1 (t))\leq -\epsilon_0 \left( c_3 +\frac{1}{c_*}\right) \left(d_1 E_1^{\prime} (t) +{\tilde d}_1 E_{1,2}^{\prime} (t)\right)
\end{equation}
\begin{equation*}
-c_3 \frac{G_{0}\,(\varepsilon_{0}E_1 (t))}{E_1 (t)} E_{1,1}^{\prime} (t). 
\end{equation*}
Similarly, multiplying \eqref{Stab2} by $\frac{G_{0}\,(\varepsilon_{0}E_2 (t))}{\varepsilon_{0}E_2 (t)}$, using \eqref{poincareine} and the definition of $E_2$ and combining with \eqref{3.23+} and \eqref{3.23+*}, we get
\begin{equation}  \label{Stab21} 
\left[\frac{2}{c_*} -\epsilon_0 \left( c_4 +\frac{1}{c_*}\right) \left(d_2 +{\tilde d}_2\right)\right] G_{0}\,(\varepsilon_{0}E_2 (t))\leq -\epsilon_0 \left( c_4 +\frac{1}{c_*}\right) \left(d_2 E_2^{\prime} (t) +{\tilde d}_2 E_{2,2}^{\prime} (t)\right)
\end{equation}
\begin{equation*}
-c_4 \frac{G_{0}\,(\varepsilon_{0}E_2 (t))}{E_2 (t)}E_{2,1}^{\prime} (t).
\end{equation*}
Because $E_1$ and $E_2$ are non-increasing and $H_0 (s):=\frac{G_0 (s)}{s}$ is non-decreasing, then $H_{0}\,(\varepsilon_{0} E_1 )$
and $H_{0}\,(\varepsilon_{0} E_2 )$ are non-increasin, and therefore
\begin{equation}
-c_3 \frac{G_{0}\,(\varepsilon_{0}E_1 (t))}{E_1 (t)} E_{1,1}^{\prime} (t)\leq -c_3 \frac{G_{0}\,(\varepsilon_{0}E_1 (0))}{E_1 (0)}
E_{1,1}^{\prime} (t) \label{Stab12} 
\end{equation}
and
\begin{equation}
-c_4 \frac{G_{0}\,(\varepsilon_{0} E_2 (t))}{E_2 (t)} E_{2,1}^{\prime} (t)\leq -c_4 \frac{G_{0}\,(\varepsilon_{0}E_1 (0))}{E_1 (0)}
E_{2,1}^{\prime} (t). \label{Stab22} 
\end{equation} 
Choosing 
\begin{equation*} 
0<\varepsilon_{0}<\frac{2}{(c_* c_3 +1)\left(d_1 +{\tilde d}_1\right)}\,\,\hbox{in case \eqref{Eq1}},\quad
0<\varepsilon_{0}<\frac{2}{(c_* c_4 +1)\left(d_2 +{\tilde d}_2\right)}\,\,\hbox{in case \eqref{Eq2}}
\end{equation*}
and exploiting \eqref{Stab11}, \eqref{Stab21}, \eqref{Stab12} and \eqref{Stab22}, we find, for some $c_7 ,\,c_8 >0$,
\begin{equation}
G_{0}\,(\varepsilon_{0}E_1 (t))\leq -c_7 \left(E_1^{\prime} (t) + E_{1,1}^{\prime} (t)+E_{1,2}^{\prime} (t)\right)\label{Stab13} 
\end{equation}
and
\begin{equation}
G_{0}\,(\varepsilon_{0}E_2 (t))\leq -c_8 \left(E_2^{\prime} (t) +E_{2,1}^{\prime} (t)+E_{2,2}^{\prime} (t)\right). \label{Stab23} 
\end{equation}
Finally, integrating \eqref{Stab13} and \eqref{Stab23} on $[0,t]$, for $t\in\mathbb{R}_+^*$, and noting that 
$G_{0}\,(\varepsilon_{0} E_1 )$ and $G_{0}\,(\varepsilon_{0} E_2 )$ are non-increasing, we arrive at
\begin{equation*}
tG_{0}\,(\varepsilon_{0} E_1 (t))\leq\int_0^t G_{0}\,(\varepsilon_{0} E_1 (s))ds\leq c_7\left(E_1 (0) +E_{1,1} (0)+E_{1,2} (0)\right):=c_9 
\end{equation*}
and
\begin{equation*}
tG_{0}\,(\varepsilon_{0} E_2 (t))\leq\int_0^t G_{0}\,(\varepsilon_{0} E_2 (s))ds\leq c_6\left(E_2 (0) +E_{2,1} (0)+E_{2,2} (0)\right):=c_{10} .
\end{equation*}
Consequently, because $G_{0}$ is inversible and non-decreasing, we deduce that
\begin{equation*}
E_1 (t)\leq \frac{1}{\varepsilon_{0}}G_{0}^{-1}\,\left(\frac{c_9}{t}\right)\quad\hbox{and}\quad E_2 (t)\leq \frac{1}{\varepsilon_{0}}G_{0}^{-1}\,\left(\frac{c_{10}}{t}\right),
\end{equation*}
which gives \eqref{203**}, for $n=1$, with 
\begin{equation*}
G_1 =G_{0}^{-1} ,\quad\alpha_{1,2} =\max\left\{c_9 ,\frac{1}{\varepsilon_{0}}\right\}\quad\hbox{and}\quad\alpha_{2,1} =\max\left\{c_{10} ,\frac{1}{\varepsilon_{0}}\right\}. 
\end{equation*}
\vskip0,1truecm
Because $D(\mathcal{A}_1^5 )\subset D(\mathcal{A}_1^4 )$ and $D(\mathcal{A}_2^3 )\subset D(\mathcal{A}_2^2 )$, then \eqref{203**} is still valid for $n=1$, $U_0\in D(\mathcal{A}_1^5 )$ in case \eqref{Eq1} and $U_0\in D(\mathcal{A}_2^3 )$ in case \eqref{Eq2}. 
\vskip0,1truecm
By induction on $n$, \eqref{203**} holds, for any $n\in \mathbb{N}^*$. Indeed, let $n\in \mathbb{N}^*$ and  
suppose that \eqref{203**} holds, for any initial data in $D(\mathcal{A}_1^ {2n+2} )$ in case \eqref{Eq1} and $D(\mathcal{A}_2^{2n} )$ in case \eqref{Eq2}. Let $U_0\in D(\mathcal{A}_1^{2(n+1)+2} )$ in case \eqref{Eq1}, $U_0\in D(\mathcal{A}_2^{2(n+1)} )$ in case \eqref{Eq2} and $U$ the corresponding solution of \eqref{222}. We have (thanks to Theorem \ref{Theorem 1.1}) 
\begin{equation*}
U_0\in D(\mathcal{A}_1^{2(n+1)+2} )\subset D(\mathcal{A}_1^{2n+2} ),\,\, U_t (0)\in D(\mathcal{A}_1^{2(n+1)+1} )\subset D(\mathcal{A}_1^{2n+2} )\,\, \hbox{and}\,\, U_{tt} (0)\in D(\mathcal{A}_1^{2n+2})\,\,\hbox{in case \eqref{Eq1}}
\end{equation*}
and
\begin{equation*}
U_0\in D(\mathcal{A}_2^{2(n+1)} )\subset D(\mathcal{A}_2^{2n} ),\quad U_t (0)\in D(\mathcal{A}_2^{2n+1} )\subset D(\mathcal{A}_2^{2n} )\quad \hbox{and}\quad U_{tt} (0)\in D(\mathcal{A}_2^{2n})\quad\hbox{in case \eqref{Eq2}},
\end{equation*}
and then \eqref{203**} holds, for $U_0$, and implies that, for some $a_{j,n} ,\, b_{j,n} >0,\,j=1,2$, 
\begin{equation}
E_{1,j} (t)\leq a_{j,n} \,G_{n}\left(\frac{a_{j,n}}{t}\right)\quad\hbox{and}\quad E_{2,j} (t)\leq b_{j,n} \,G_{n}\left(\frac{b_{j,n}}{t}\right). \label{3.38} 
\end{equation} 
By integrating \eqref{Stab13} and \eqref{Stab23} over $[T,2T]$, for $T>0$, noting that $G_{0}(\epsilon_0 E_j )$ is non-increasing and using \eqref{3.38}, we get, for some $d_{j,n} >0,\,j=1,2$,
\begin{equation}
TG_{0}(\epsilon_0 E_1 (2T)) \leq \displaystyle\int_{T}^{2T}G_0 (\epsilon_0 E_1 (t))\,dt \leq c_7 (E_1 (T)+E_{1,1} (T)+ E_{1,2} (T)) 
\leq d_{1,n} G_{n}\left(\frac{d_{1,n}}{T}\right)\label{3.39} 
\end{equation}
and
\begin{equation}
TG_{0}(\epsilon_0 E_2 (2T)) \leq \displaystyle\int_{T}^{2T}G_0 (\epsilon_0 E_2 (t))\,dt\leq
c_8 (E_2 (T)+E_{2,1} (T)+ E_{2,2} (T)) \leq d_{2,n} G_{n}\left(\frac{d_{2,n}}{T}\right).\label{3.39+} 
\end{equation}
Therefore, since $G_0$ is non-decreasing, 
\begin{equation*}
E_1 (2T)\leq\frac{1}{\varepsilon_{0}}G_{0}^{-1}\,\left(\frac{2d_{1,n}}{2T}G_{n}\left(\frac{2d_{1,n}}{2T}\right)\right)\quad\hbox{and}\quad E_2 (2T)\leq \frac{1}{\varepsilon_{0}}G_{0}^{-1}\,\left(\frac{2d_{2,n}}{2T}G_{n}\left(\frac{2d_{2,n}}{2T}\right)\right),
\end{equation*}
that is 
\begin{equation*}
E_1 (t)\leq\alpha_{1,n+1} G_{n+1}\left(\frac{\alpha_{1,n+1}}{t}\right)\quad\hbox{and}\quad E_2 (t)\leq \alpha_{2,n+1} G_{n+1}\left(\frac{\alpha_{2,n+1}}{t}\right),\quad t>0,
\end{equation*}
where 
\begin{equation*}
G_{n+1} (s)=G_1 (sG_n (s)), \quad\alpha_{1,n+1}=\max\left\{\frac{1}{\varepsilon_{0}},2d_{1,n}\right\}\quad 
\hbox{and}\quad\alpha_{2,n+1} =\max\left\{\frac{1}{\varepsilon_{0}},2d_{2,n}\right\},
\end{equation*}
which leads to \eqref{203**}, for $n+1$ instead of $n$. This ends the proof of Theorem \ref{theorem1}. 

\section{General comments}

1. If $g_0 =0$, then $g\equiv 0$, and therefore, \eqref{St2} and \eqref{St3} lead to $E_1 (t)=E_1 (0)$ and $E_2 (t)=E_2 (0)$, for all
$t\in \mathbb{R}_+$. So the presence of the memory term is necessary to get the stability of \eqref{222}. 
\vskip0,1truecm
2. One can consider the more general form of the first equations in \eqref{Eq1} and \eqref{Eq2} by considering
\begin{equation}\label{Ge1}
iy_{t} (x,t)+ A y(x,t)- \displaystyle\int_0^{\infty}\,g (s) B y (x,t-s) \,ds = 0
\end{equation}
and
\begin{equation}\label{Ge2}
iy_{t} (x,t)+ A y(x,t) +c(x)\displaystyle\int_0^{\infty}\,g (s) y (x,t-s) \,ds = 0,
\end{equation}
where $c$ is non-negative real valued functions belonging to $C ({\bar{\Omega}})$,
\begin{equation*}
A=\sum_{k,j=1}^d \frac{\partial}{\partial x_k} \left(a_{kj}\frac{\partial}{\partial x_j}\right) \quad\hbox{and}\quad B=\sum_{k,j=1}^d 
\frac{\partial}{\partial x_k} \left(b_{kj}\frac{\partial}{\partial x_j}\right) 
\end{equation*}
such that $a_{kj}$ and $b_{kj}$ are real valued functions belonging to $C^1 ({\bar{\Omega}})$ satisfying, for some $a_0 ,\,b_0 ,\,c_0 >0$, 
\begin{equation*}
a_{kj} (x)=a_{jk} (x),\quad \sum_{k,j=1}^d a_{kj} (x)\epsilon_k \epsilon_j \geq a_0 \sum_{k=1}^d \epsilon_k^2 ,
\end{equation*}
\begin{equation*}
b_{kj} (x)=b_{jk} (x),\quad\sum_{k,j=1}^d b_{kj} (x)\epsilon_k \epsilon_j \geq b_0 \sum_{k=1}^d \epsilon_k^2  
\end{equation*}
and $c(x)\geq c_0$, for any $(\epsilon_1,\cdots , \epsilon_d ) \in  \mathbb{R}^d$ and $x\in \Omega$. An abstract form including \eqref{Ge1} and \eqref{Ge2} can be also considered by taking \eqref{Ge1} with self-adjoint linear positive definite operators $A:\,D(A)\to H$ and $B:\,D(B)\to H$ and a Hilbert space $H$ with dense and compact embeddings $D(A)\subset D(B)\subset H$ such that there exist positive constants $a_0$ and $b_0$ satisfying
\begin{equation*}
b_0 \Vert v\Vert^2 \leq \Vert B^{\frac{1}{2}} v\Vert^2 \leq a_0\Vert A^{\frac{1}{2}} v\Vert^2,\quad\forall v\in D(A^{\frac{1}{2}}). 
\end{equation*}
\vskip0,1truecm
3. Our results hold if we consider a domaine $\Omega$ not necessarily bounded but of a finite measure;
so Poincar\'e's inequality \eqref{poincareine} is still applicable. Howover, considering $\mathbb{R}^d$ or a domaine with infinite measure is a nice open question. 
\vskip0,1truecm
4. We can add to \eqref{Eq1} and \eqref{Eq2} a linear term of the form $b(x)y(x,t)$, where $b$ is a real valued function belonging to 
$C({\bar{\Omega}})$ and satisfying $\Vert b\Vert_{\infty}<\frac{a}{c_*}$ ($c_*$ is the Poincar\'e's constant defined in \eqref{poincareine}). It will be nice to study the case where a non-linear term of the form $h(\vert y(x,t)\vert) y(x,t)$ is added to the first equations in \eqref{Eq1} and \eqref{Eq2}, where $h:\,\mathbb{R}_+\to \mathbb{R}$ is a given function. The difficulty in the non-linear case is that \eqref{St2*} is not satisfied.
\vskip0,1truecm
5. Our hypothesis \eqref{203} allows $g$ to have a decay rate at infinity faster than $\frac{1}{s^3}$. The case of $g$ having a decay rate at infinity between $\frac{1}{s^3}$ and $\frac{1}{s}$ is open. In the case of hyperbolic systems considered in \cite{gues1} and \cite{gues2}, it was assumed that (instead of \eqref{203}) 
\begin{equation*}
\displaystyle\int_{0}^{\infty}\frac{g (s)}{G^{-1}\left(-\,g' (s)\right)}\ ds +\sup_{s\in\mathbb{R}_{+}}\frac{g (s)}{G^{-1}
\left(-\,g' (s)\right)} < \infty ,
\end{equation*}
which allows $g$ to have a decay rate at infinity arbitrarily close to $\frac{1}{s}$, and the obtained decay rate for the corresponding energy $E$ was better, more precisely, it was proved in \cite{gues1} and \cite{gues2} that (instead of \eqref{203**}), 
for some $\alpha_{n} >0$,
\begin{equation*}
E(t)\leq \alpha_{n} \,G_{n}\left(\frac{\alpha_{n}}{t}\right),
\end{equation*}
for any $n\in \mathbb{N}^*$ and any $U_0\in D(\mathcal{A}^n)$. This is because in the case of hyperbolic systems, the adequate variable $\eta^t$ using to treat the infinit memory is defined by $\eta^t (x,s)=y(t)-y(t-s)$ (instead of \eqref{etaname}). 

\section{Numerical Examples}

In order to numerically corroborate the asymptotic behavior and exponential decay of energy, we will show some examples in dimension 1, that is, the domain is simply an interval $ (0, L) $.

Since the problems \eqref{Eq1} and \eqref{Eq2} are linear, it is convenient to directly use Fourier series to approximate problem (1.1), that is, we use the method of separation of variables to assume that the solutions are shared as:
\begin{equation}
\label{SinusSerie}
y(x,t)=\sum_{k=1}^\infty B_k(t)\sin(2k\pi x),
\end{equation}
given the Dirichlet initial conditions.
In this case, the constants $ B_k $, verify the following differential-integral equations:
\begin{align}\label{EDB1}
\begin{cases}
i B_k'- a \dfrac{4\pi^2k^2}{L^2} B_k+ i\displaystyle\int_0^{\infty}\,f (s) \dfrac{4\pi^2k^2}{L^2}  B_k (t-s) \,ds = 0, \quad &x\in \Omega,\,t\in \mathbb{R}_+^* :=(0,\infty),\\
B_k(-t)=B_k^0(t),\quad &t\in \mathbb{R}_+ :=[0,\infty)
\end{cases}
\end{align}
and 
\begin{align}\label{EDB2}
\begin{cases}
i B_k'- a \dfrac{4\pi^2k^2}{L^2} B_k+ i\displaystyle\int_0^{\infty}\,f (s)  B_k (t-s) \,ds = 0, \quad &x\in \Omega,\,t\in \mathbb{R}_+^* ,\\
B_k(-t)=B_k^0(t),\quad &t\in \mathbb{R}_+ :=[0,\infty)
\end{cases}
\end{align}
respective approximations of \eqref{Eq1} and \eqref{Eq2}.
In this case $ B_k^0 $ represent the coefficients of the series of sines \eqref{SinusSerie} for the initial condition function $y_0 (x,t)$.
In order to solve these integral differential equations, we will use Heun's method, which, in addition to being a second order scheme, preserves energy for the linear Schrodinger equation without dissipative term. This is the Crank-Nicolson version for this single variable equation:
\begin{align}\label{scheme1}
\begin{cases}
i \dfrac{B_k^{n+1}-B_k^n}{\delta t}- a \dfrac{4\pi^2k^2}{L^2} B^{n+\frac{1}{2}}_k+ i\displaystyle\dfrac{4\pi^2k^2}{L^2} \sum_{m=0}^N\,\delta t f^m  B^{n-m+\frac{1}{2}}_k \,ds = 0, \\
B_k^{-n}=B_k^{0,n},\quad k=1,\ldots,K, \quad n=0,\ldots,N,
\end{cases}
\end{align}
and 
\begin{align}\label{scheme2}
\begin{cases}
i \dfrac{B_k^{n+1}-B_k^n}{\delta t}- a \dfrac{4\pi^2k^2}{L^2} B^{n+\frac{1}{2}}_k+ i\displaystyle \sum_{m=0}^N\,\delta t f^m  B^{n-m+\frac{1}{2}}_k \,ds = 0,  \\
B_k^{-n}=B_k^{0,n},\quad k=1,\ldots,K, \quad n=0,\ldots,N,
\end{cases}
\end{align}
where $B^{n+\frac{1}{2}}_k=\dfrac{B^{n+1}_k+B^{n}_k}{2}$.
\begin{figure}[hbtp!]
\includegraphics[width=7cm]{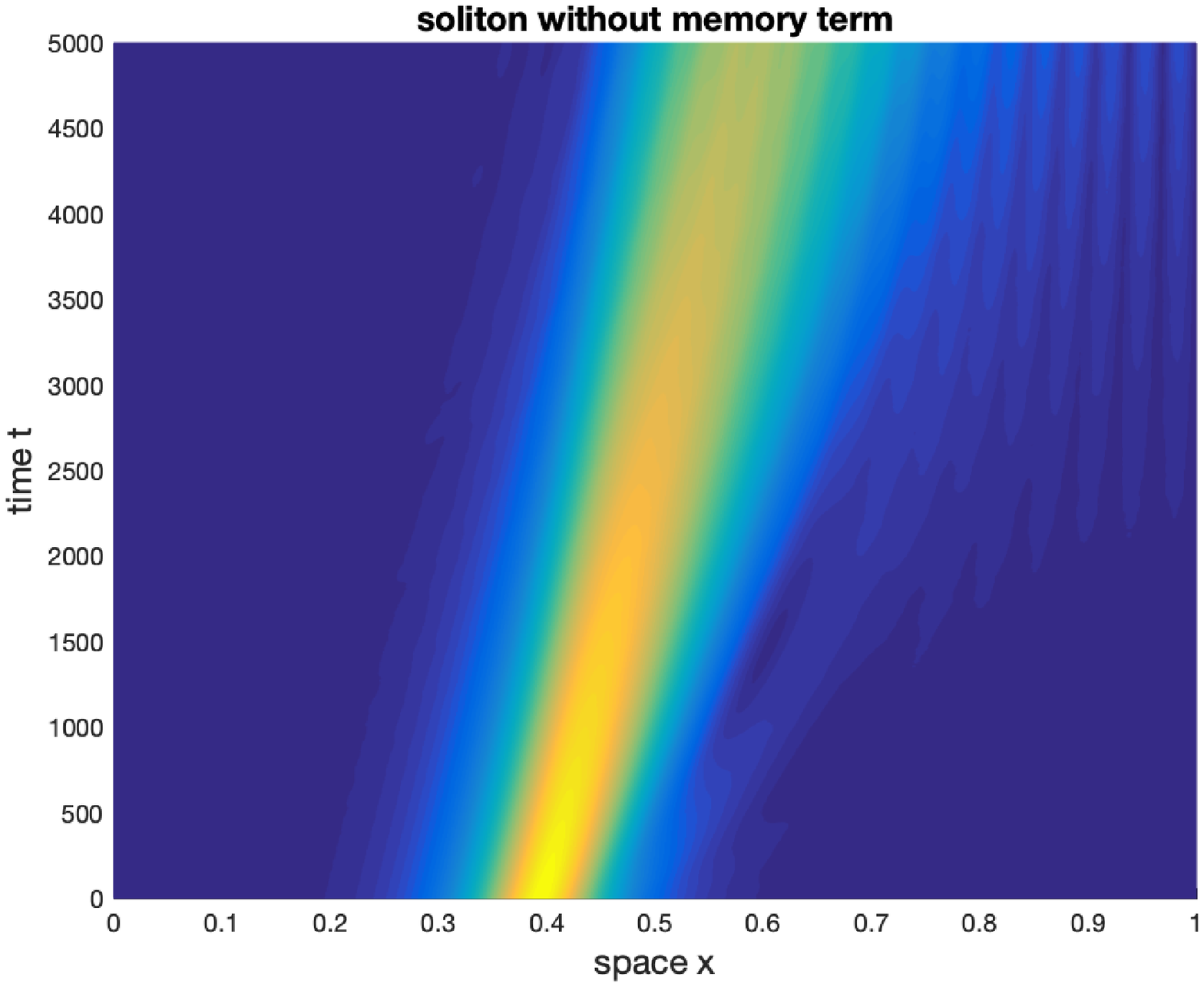}
\includegraphics[width=7cm]{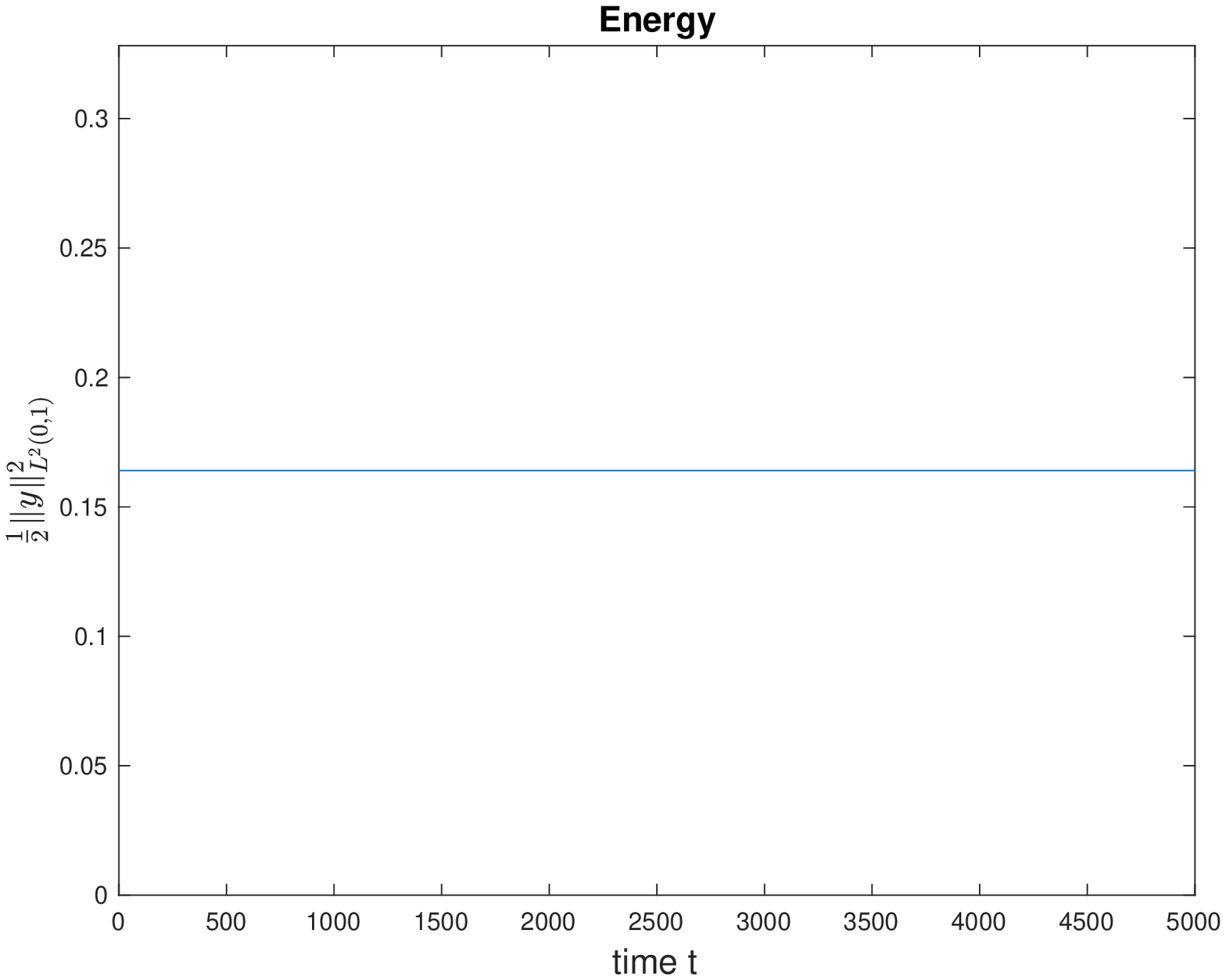}
\caption{Numerical solution of the Schr\"odinger equation without memory term (graph on the left) and Energy (on the right).}
\end{figure}

\subsection{Example}
We consider a  initial condition, constant with respect to time:
$$
y_0(x,t)= \dfrac{A e^{i \lambda x}}{\cosh\left(\dfrac{x-x_1}{x_0}\right)},
\qquad x\in(0,L), \ t>0,
$$
with $A=4$, $\lambda=7$, $x_0=\dfrac{1}{2A\sqrt{\lambda}}$,
$x_1=0.4$, $L=1$.
We assume $a=1$, and $f(t)=-\displaystyle \int^t_0 g_1(s)\,ds$ or $f(t)=-\displaystyle\int^t_0 g_2(s)\,ds$,
given by formula \eqref{g2} and \eqref{g1}, respectively.
We choose $d_1=d_2=10.000$, $q_1=1$, and $q_2=4$.
The simulations are done with $K = 2^{10}$ and $N = 20000$ ($T = 2000$, $\delta t = 0.1$).

In Figure 1, the solution in space and time of the linear Schrodinger equation is observed. Given the initial condition, this should be a soliton for the nonlinear Schrodinger on the entire real line \cite{hillion}. However, given the linear equation and the bounded domain with Dirichlet edge conditions, the soliton dissolves by interacting with part of the scattered and reflected waves at the edge. The important thing is that in this case the energy is completely conserved as shown in Figure 2 (blue line in the graph on the left).

\begin{figure}[hbtp!]
\includegraphics[width=7cm]{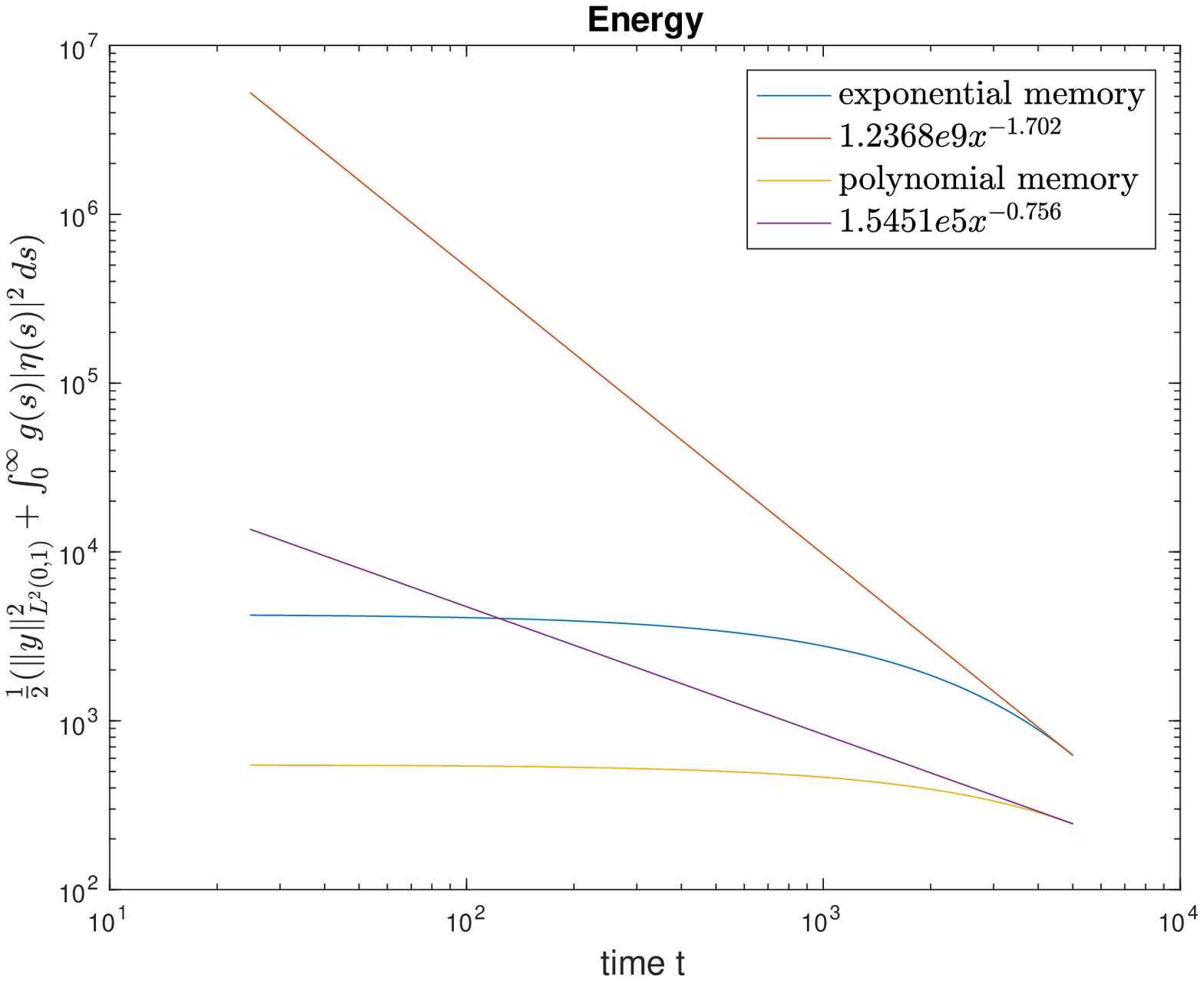} 
\includegraphics[width=7cm]{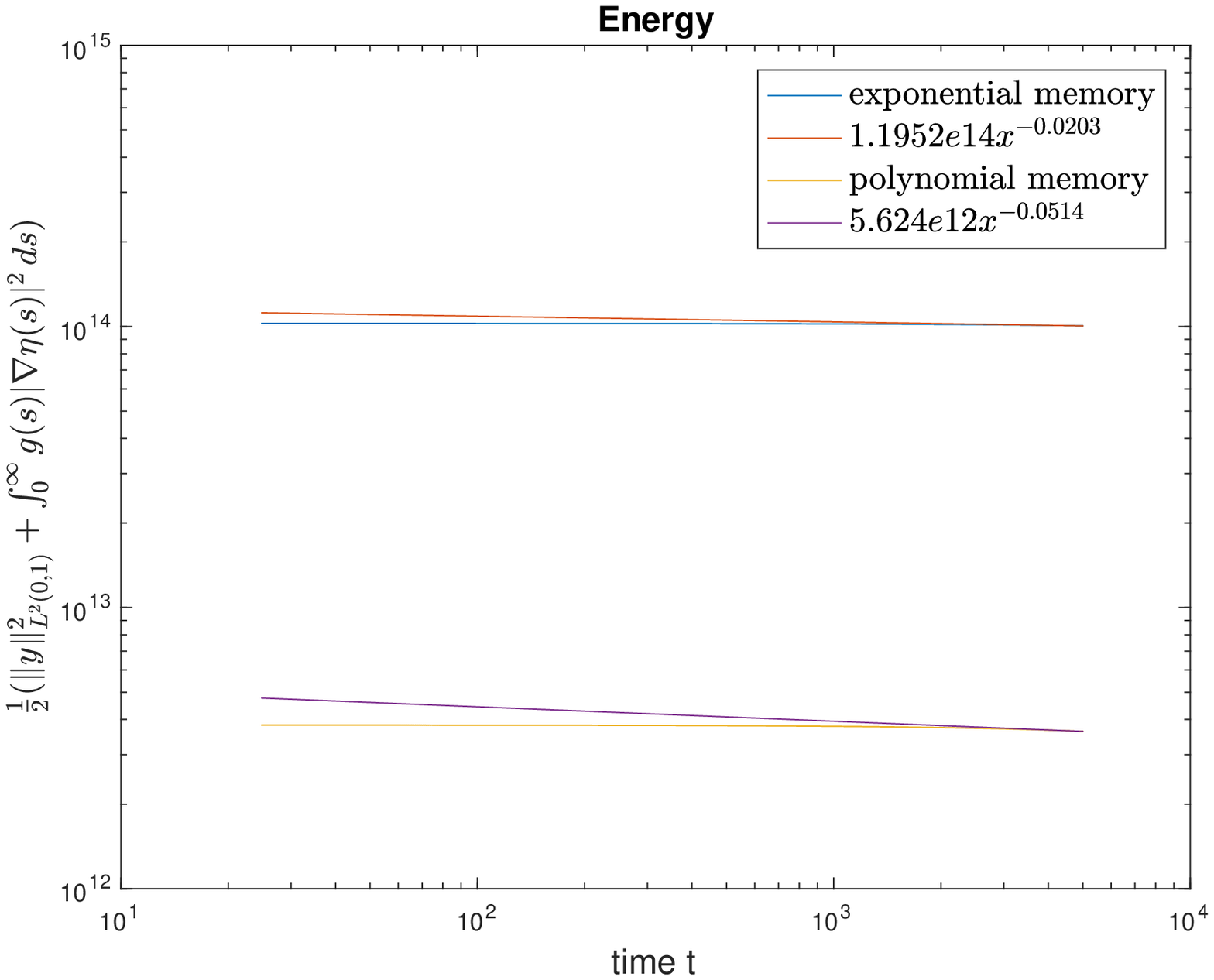}
\caption{Exponential decay of the Energy. On the left: logarithmic scale graphics of the energies
for the  Schr\"odinger equation with memory term 
given by equation \eqref{Eq2}. On the right: logarithmic scale graphics of the energies
for the equation \eqref{Eq1}.}
\end{figure}

Using the Parseval Identity, the energies \eqref{En1} and \eqref{En2} are naturally discretized by 
\begin{equation}\label{EnD1}
E_{1,\delta}^n= \frac{L}{4}  \left(\sum_{k=1,\ldots,K}\left(B_k^n\right)^2 
+\displaystyle \sum_{k=1}^K\sum_{m=1}^N \delta t \dfrac{4\pi^2k^2}{L^2} g^m 
\left(\eta_k^{m,n}\right)^2\right),
\end{equation}
and
\begin{equation}\label{EnD2}
E_{2,\delta}^n =
\frac{L}{4}  \left(\sum_{k=1,\ldots,K}\left(B_k^n\right)^2 
+\displaystyle \sum_{k=1}^K\sum_{m=1}^N \delta t g^m 
\left(\eta_k^{m,n}\right)^2\right),
\end{equation} 
where $\displaystyle \eta_k^{m,n}=\sum_{\ell=n-m}^n \delta t B_k^\ell$.
In Figure 2, graphs of the discrete energies  \eqref{EnD1} and \eqref{EnD2} are observed in  logarithmic scale, to compare the exponential decay with the different types of memory. The case without memory appears in the graph on the left in blue. It is observed that the exponential function $ g_1 $ causes the decay with the highest rate for the memory term equal to $\displaystyle\int_0^\infty \frac{d_1}{q_1} e^{-q_1 s}y(t-s)\,ds$. The graphs on the right correspond to the memory terms equal to $-\displaystyle\int_0^\infty f(s)\Delta y(t-s)\,ds$, and in addition to having greater energy, it is observed that they decay exponentially more slowly with a smaller rate.

\begin{figure}[hbtp!]
\begin{subfigure}[b]{.4\textwidth}
\includegraphics[width=7cm]{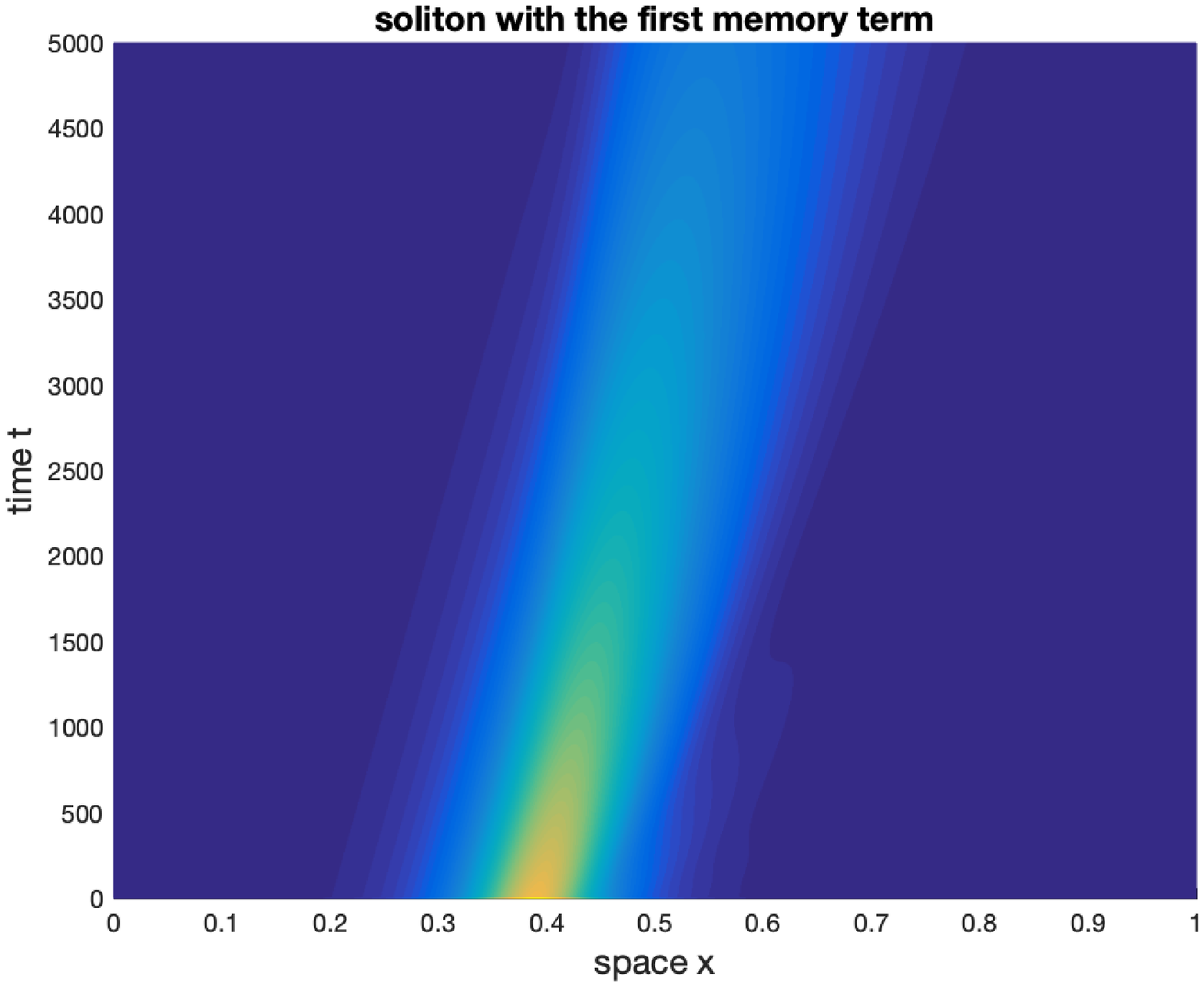} 
\caption{$\displaystyle\int_0^\infty \frac{d_1}{q_1} e^{-q_1 s}y(t-s)\,ds$}
\end{subfigure}
\begin{subfigure}[b]{.4\textwidth}
\includegraphics[width=7cm]{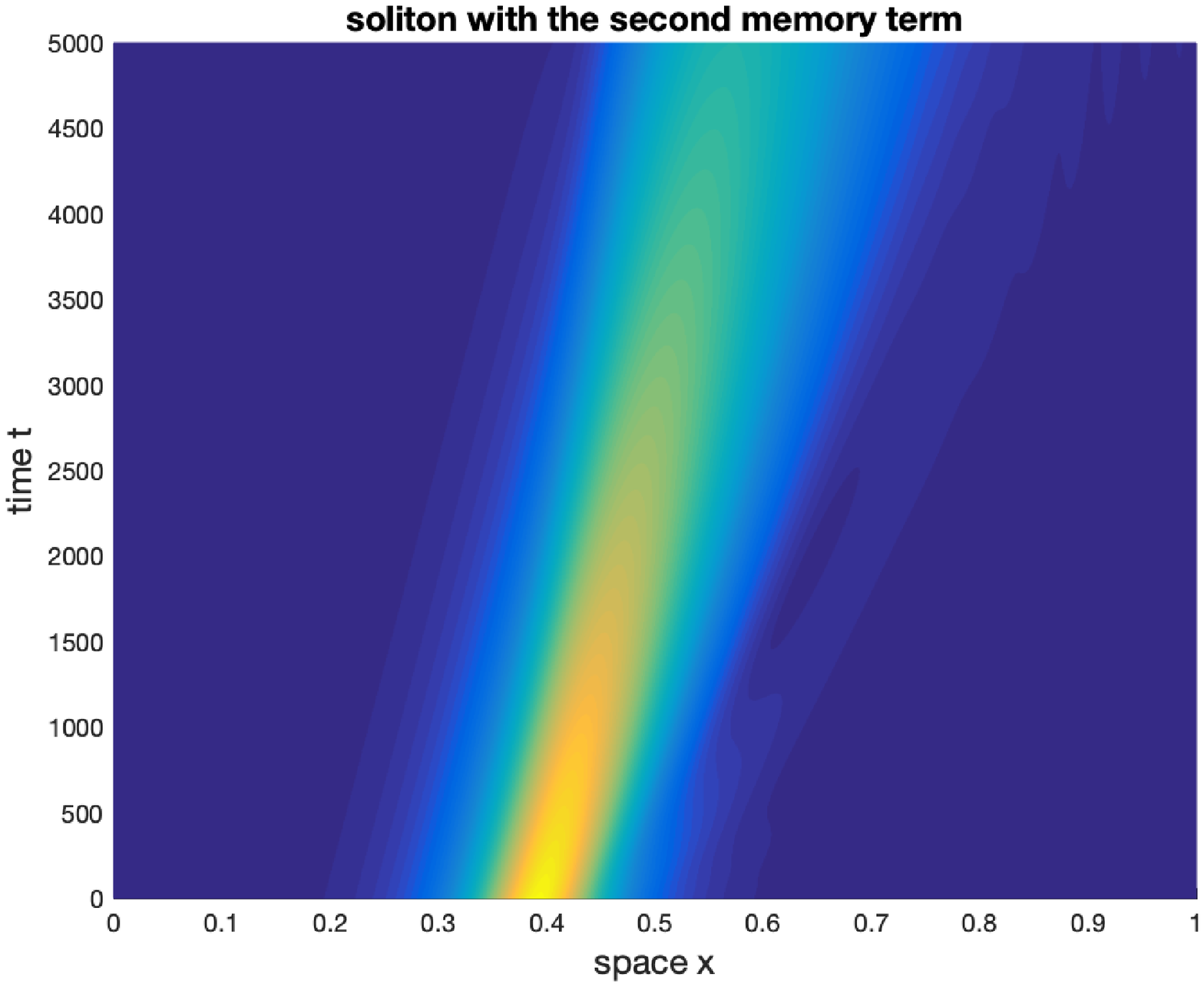} 
\caption{$\displaystyle\int_0^\infty  \dfrac{d_2/(q_2-1)}{(1+s)^{q_2-1}}y(t-s)\,ds$}
\end{subfigure}
\\
\begin{subfigure}[b]{.4\textwidth}
\includegraphics[width=7cm]{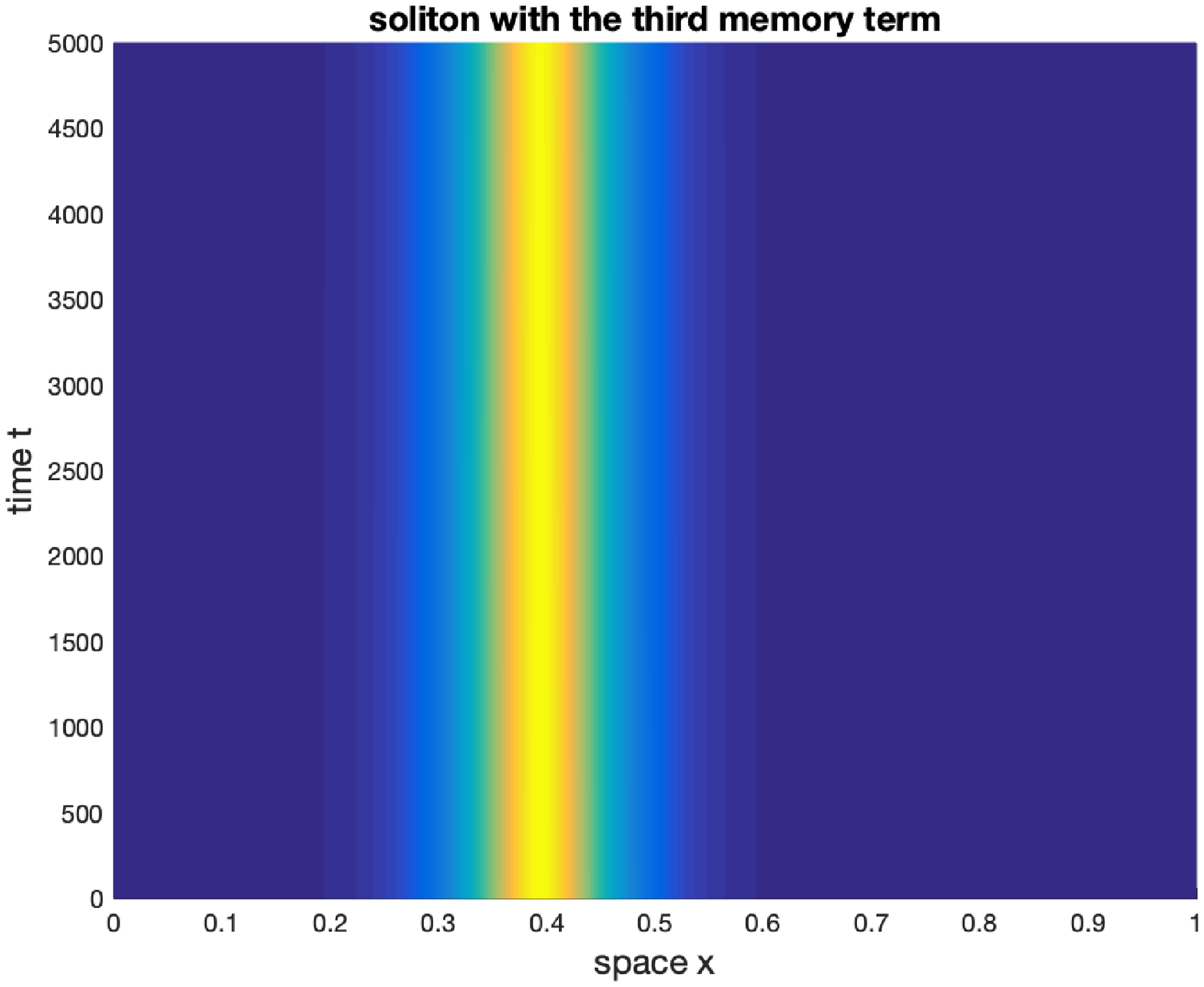}
\caption{$-\displaystyle\int_0^\infty \frac{d_1}{q_1}  e^{-q_1 s}\Delta y(t-s)\,ds$}
\end{subfigure}
\begin{subfigure}[b]{.4\textwidth}
\includegraphics[width=7cm]{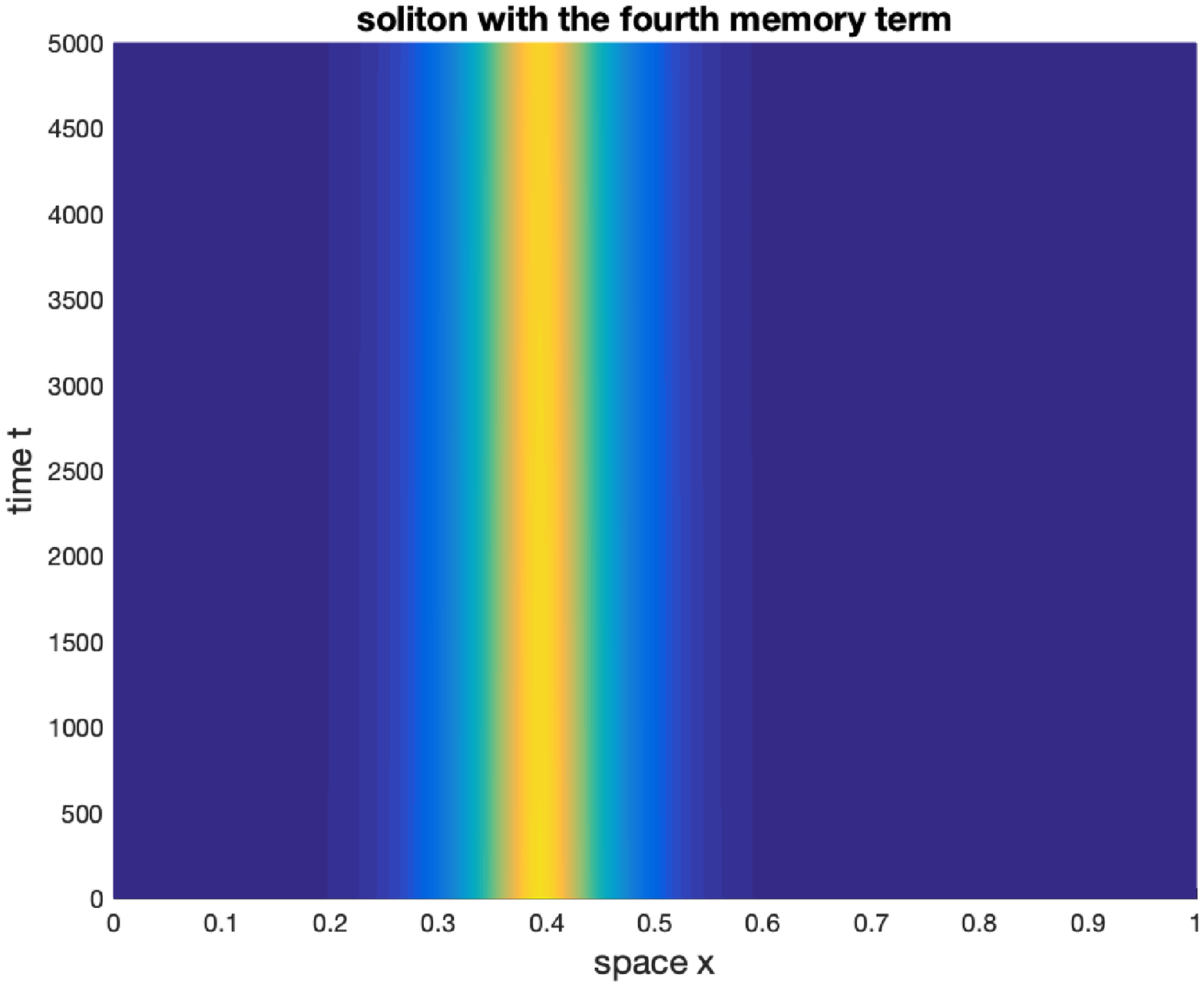}
\caption{$\displaystyle-\int_0^\infty  \dfrac{d_2/(q_2-1)}{(1+s)^{q_2-1}}\Delta y(t-s)\,ds$}
\end{subfigure}
\caption{Solution of the Schr\"odinger equation with different memory terms.}
\end{figure}

Finally, in Figure 3, the asymptotic behavior of the different solutions is observed, with the 4 different types of memory. In the case of memory terms of the form $\displaystyle\int_0^\infty f(s) y(t-s)\,ds$, there is a fading of the soliton, much more marked than that of the Schrodinger equation without memory (graphs on the left), on the other hand when memory terms of the form $-\displaystyle\int_0^\infty f(s)\Delta y(t-s)\,ds$ are considered , the soliton remains visually intact, however, Figure 2 indicates that the energy clearly decays
(graphs on the right).

\vskip0,2truecm
{\bf Acknowledgment.} This work was initiated during the visit in July-August 2017 of the third author to Concepci\'{o}n University, Chile, and finished during the visit of the fourth author in June 2018 to Lorraine - Metz university, France, and the visits in August 2018 of the third author to Concepci\'{o}n university, Chile, and Maringa University, Brazil. The authors thank Concepci\'{o}n, Lorraine - Metz and Maringa Universities for their kind support and hospitality.

\end{document}